\documentclass[12pt,leqno]{article}
\usepackage{amsmath}
\usepackage{amssymb}
\usepackage{amscd}
\usepackage{amsfonts}
\usepackage{amsthm}
\usepackage[matrix,arrow,curve]{xy}

\newtheorem{df}{Definition}
\newtheorem{theorem}{Theorem}
\newtheorem{prop}[theorem]{Proposition}
\newtheorem{lemma}[theorem]{Lemma}
\newtheorem{conj}[theorem]{Conjecture}
\newtheorem{cor}[theorem]{Corollary}

\newtheorem{rem}{Remark}

\newcommand{\Dj}{\hbox to 8pt{\raisebox{.4\height}{-}\hss D}}
\newcommand{\dJ}{\,{\hbox to 4pt{\raisebox{.75\height}{-}\hss d}}}
\newcommand{\id}{\ensuremath{{\rm id}}}

\newcommand{\Hom}{\ensuremath{\mathop{\rm Hom}\nolimits}}
\newcommand{\im}{\ensuremath{\mathop{\rm Im}\nolimits}}

\newcommand{\ac}{\ensuremath{\mathcal A}}
\newcommand{\mc}{\ensuremath{\mathcal M}}

\newcommand{\hc}{\ensuremath{\mathcal H}}

\newcommand{\nc}{\ensuremath{\mathcal N}}
\newcommand{\rc}{\ensuremath{\mathcal R}}

\newcommand{\xc}{\ensuremath{\mathcal X}}
\newcommand{\yc}{\ensuremath{\mathcal Y}}

\newcommand{\eqdef}{\stackrel{\rm def}{=}}
\newcommand{\lar}{\longrightarrow}
\newcommand{\lal}{\longleftarrow}

\begin{document}
\begin{title}
{\vspace{-3cm}
\begin{flushright}
\normalsize Preprint: ITEP-TH-108/05
\end{flushright}
\vspace{5cm}
Pairings in Hopf-cyclic cohomology of algebras and coalgebras with coefficients.}
\end{title}
\author{I.~Nikonov \thanks{dpt. of Differential Geometry, MSU, Leninskie Gory, Moscow, Russia}, G.~Sharygin\thanks{ITEP, ul. B. Cheremushkinskaja, 25, Moscow, Russia}\ \thanks{The second author was partly supported by the grant RFFI 04-01-00702}}
\date{}
\maketitle

{\begin{abstract}
This paper is concerned with the theory of cup-products in Hopf-type cyclic cohomology of algebras and coalgebras.
Here we give detailed proofs of the statements, announced in our
previous paper \cite{nashe}. We show that the cyclic cohomology of a
coalgebra can be obtained from a construction involving noncommutative Weil
algebra. Then we use a generalization of Quillen and Crainic's construction
(see \cite{Quil} and \cite{Crain}) to define the cup-product. We discuss the
relation of the introduced cup-product and $S$\/-operations on cyclic
cohomology. After this we describe the relation of this type of product and
bivariant cyclic cohomology. In the last section we briefly discuss the relation of our constructions with that of \cite{KhalRan}.
\end{abstract}}
\bigskip

\tableofcontents

\section{Definitions and notation}
Throughout the paper \hc\ will denote a fixed Hopf algebra with
invertible antypode over a fixed characteristic zero field $\Bbbk$
and \mc\ -- a stable anti-Yetter-Drinfeld module over \hc. We assume
that \mc\ is a left \hc\/-comodule and right module. Throughout the
text we shall use the standard (Sweedler's) notation with
superscripts for all the comultiplications and coactions, e.g.
$\Delta(h)=h^{(1)}\otimes h^{(2)}$ for all $h\in\hc$ and
$\Delta_\mc(m)=m^{(-1)}\otimes m^{(0)}$ for every $m\in\mc$. Under
this assumptions the anti-Yetter-Drinfeld condition takes form
\begin{equation}\label{aYDc}
(mh)^{(-1)}\otimes(mh)^{(0)}=S(h^{(3)})m^{(-1)}h^{(1)}\otimes
m^{(0)}h^{(2)},
\end{equation}
and the stability of \mc\ means that $m^{(0)}m^{(-1)}=m$ for all $m\in\mc$.

Now recall the definition of the algebras and coalgebras cyclic
cohomology with coefficients in stable anti Yetter-Drinfeld modules.

Let $C$ be a coalgebra over the same field. Suppose \hc\/-acts on $C$ from the left
$\hc\otimes C\to C$, in a way that respects the coalgebraic structure of $C$, i.e.
$(hc)^{(1)}\otimes(hc)^{(2)}=h^{(1)}c^{(1)}\otimes h^{(2)}c^{(2)}$ for all
$h\in\hc,\ c\in C$.

Recall the definitions of the Hopf-type cyclic cohomology of $C$ with
coefficients in \mc, $HC^\ast_\hc(C,\,\mc)$, given in \cite{HKRS}.

First of all, one considers paracocyclic module $C^\ast(C,\,\mc)$:
\begin{align}
&C^n(C,\,\mc)=\mc\otimes C^{\otimes n+1},\\
\intertext{and the cyclic operations are defined by the formulas}
\delta_i(m\otimes c_0\otimes\dots\otimes c_{n})&=\begin{cases}
   m\otimes c_0\otimes\dots\otimes c_i^{(1)}\otimes
c_i^{(2)}\otimes\dots\otimes c_n,&\ 0\le i\le n,\\
   m^{(0)}\otimes c_0^{(1)}\otimes c_1\otimes\dots\otimes c_n\otimes
m^{(-1)}c_0^{(2)},&\ i=n+1,
                                                   \end{cases}\\
\sigma_i(m\otimes&c_0\otimes\dots\otimes c_{n})=m\otimes c_0\otimes\dots\otimes\epsilon(c_i)\otimes\dots\otimes
c_{n},\\
\tau_n(m\otimes&c_0\otimes\dots\otimes c_{n})=m^{(0)}\otimes c_1\otimes\dots\otimes c_n\otimes
m^{(-1)}c_0.
\end{align}
Recall that ``paracocyclic" means that all the usual cocyclic
relations are satisfied, probably except for $\tau_n^{n+1}=1$. $C$
being a \hc\/-module, one can extend the action of the Hopf algebra
to the tensor power of $C$ diagonally and consider the factor-module
$C^\ast_\hc(C,\,\mc),\ C^n_\hc(C,\/\mc)=\mc\otimes_\hc (C^{\otimes
n+1})$. Now it is easy to show that the paracocyclic operations,
introduced above, can be pulled down to $C^\ast_\hc(C,\,\mc)$ iff
\mc\ is anti-Yetter-Drinfeld. And if \mc\ is also stable, then
$C^\ast_\hc(C,\,\mc)$ with the operations restricted on it from
$C^\ast(C,\,\mc)$ is cocyclic.

By definition Hopf-type cyclic (respectively, periodic) cohomology
of $C$ with coefficients in \mc, $HC^\ast_\hc(C,\,\mc)$ (resp.
$HP^\ast_\hc(C,\,\mc)$) is the cyclic (resp. periodic) cohomology of
the cocyclic module $C^\ast_\hc(C,\,\mc)$. This means that one
introduces the mixed complex with differentials $b$ and $B$,
associated in a usual way with the cocyclic module
$C^\ast_\hc(C,\,\mc)$ (see, e.g. \cite{Loday}) and takes the
cohomology of the corresponding total complex (resp., periodic
super-complex).

Similarly, let $A$ be a (left) Hopf-module algebra over \hc, i.e. there's an action
$\hc\otimes A\to A$, such that for all $a,\ b\in A$ and all $h\in\hc$ one has
$h(ab)=h^{(1)}(a)h^{(2)}(b)$. The following construction is also taken from
\cite{HKRS}.

Consider the (paraco)cyclic module $C^\ast(A,\,\mc)$:
\begin{align}
C^n&(A,\,\mc)=\mathrm{Hom}(\mc\otimes A^{\otimes n+1},\,k)\\
\intertext{where $\mathrm{Hom}(A,\,B)$ is the space of $k$\/-linear
homomorphisms from $A$ to $B$. The (paraco)cyclic operations on
$C^n(A,\,\mc)$ are given by the following formulae (where for the
sake of brevity we have substituted commas for the tensor product
signs} \delta_if(m,a_0,\dots,a_{n+1})&=\begin{cases}
                                  f(m,a_0,\dots,a_ia_{i+1},\dots,a_{n+1}),\ &0\le{i}\le{n},\\
                                  f(m^{(0)},S^{-1}(m^{(-1)})(a_{n+1})a_0,\dots,a_n),\ &i=n+1,
                                \end{cases}\\
\sigma_if(m,a_0,\dots,a_{n-1})&=f(m,a_0,\dots,1,\dots,a_{n-1}),\ 0\le i\le n,\\
\intertext{where $1$ stands in the $i$\/-th place and}
\tau_nf(m,a_0,\dots,a_n)&=f(m^{(0)},S^{-1}(m^{(-1)})(a_{n}),a_0,\dots,a_{n-1}).
\end{align}
Again, all the cyclic relations except $\tau_n^{n+1}=1$ are satisfied. One
passes from {\em para\/}-cocyclic to cocyclic module by taking space \hc\/-linear homomorphisms $\mathrm{Hom}_\hc(\mc\otimes A^{\otimes
n+1},\,k)$, where we let \hc\ act on $k$ on the left via the counit and the
left \hc\/-action on $\mc\otimes A^{\otimes n+1}$ is given by
$$
h\cdot(m\otimes a_0\otimes\dots\otimes a_n)=mS(h^{(1)})\otimes
h^{(2)}a_0\otimes\dots\otimes h^{(n+2)}a_n.
$$
The (co)cyclic module which is obtained in this way is denoted
$C^*_\hc(A,\,\mc)$ and its cyclic (resp. Hochschild, resp. periodic)
cohomology are called the Hopf-type cyclic (resp. Hochschild, resp. periodic)
cohomology of $A$ with coefficients in \mc. They are denoted by
$HC^*_\hc(A,\,\mc)$ (resp. $H^*_\hc(A,\,\mc)$, resp. $HP^*_\hc(A,\,\mc)$).

Suppose now that $C$ acts on $A$ so that
\begin{align}
\label{eqch1}
c(ab)&=c^{(1)}(a)c^{(2)}(b),\\
\intertext{and}
\label{eqch2}
(h(c))(a)&=h(c(a))
\end{align}
for all $a,b\in A,\ c\in C,\ h\in\hc$. In paper \cite{KhalRan} there was defined a pairing
\begin{equation}
\label{eqf}
HC^p_\hc(C,\,\mc)\otimes HC^q_\hc(A,\,\mc)\stackrel{\sharp}{\longrightarrow}
HC^{p+q}(A),
\end{equation}
extending to higher dimensions the Connes-Moscovici characteristic map
(see \cite{CM1},\cite{CM2} and \cite{HKRS}) $\gamma:HC_\hc^*(C,\,\mc)\to
HC^*(A)$, constructed for an equivariant trace $\gamma:\mc\otimes A\to k$ (in Connes' and
Moscovici's papers this map was defined olny for the $1$\/-dimesional modules
$\mc={}_\sigma\!k^\delta$ for a modular pair in involution $(\sigma,\,\delta)$).

On the other hand when $\mc={}_\sigma\!k^\delta$ the map $\gamma$ was generalized
to higher equivariant traces by
Crainic in \cite{Crain}. Methods used in that paper are quite different from
those of \cite{KhalRan}. In the paper \cite{nashe} the authors gave a brief
outline of the recipe which allows one extend the methods used by Crainic to
obtain another construction of the pairing, similar to \eqref{eqf}. We shall denote this pairing by $\sharp'$.

The present paper is
devoted to, first of all, giving a detailed proofs of the statements only
formulated and/or sketched in \cite{nashe}. In particular, we discuss the $S$\/-operation
relations, verified by the map $\sharp'$. Second, we propose another
construction of a pairing, similar to $\sharp$, this time follwoing the
methods of the book \cite{Loday}. 

\section{Non-commutative Weil algebra and its cohomologies}
Let $C$ be a coalgebra over a characteristic zero field $k$. Recall
the definition of the Weil algebra of $C$ given at \cite{nashe},
which is a straightforward generalization of the definition in
\cite{Crain}, where it is given for Hopf algebras.
\begin{df}
One calls the ``Weil algebra of a coalgebra $C$" the free
differential graded algebra (without unit) generated by elements
$i_c,\ {\rm deg}\,i_c=1$ and $w_c,\ {\rm deg}\,w_c=2$, where both
symbols are linear in $c\in C$. The differential is given by
\begin{align}
\partial i_c&=w_c-i_{c^{(1)}}i_{c^{(2)}}\\
\partial w_c&=w_{c^{(1)}}i_{c^{(2)}}-i_{c^{(1)}}w_{c^{(2)}}.
\end{align}
Weil algebra of a coalgebra $C$ is denoted by $W(C)$. And the symbol
$I(C)$ is reserved for the canonical ideal in $W(C)$ generated by
the elements $w_c,\ c\in C$.
\end{df}
For better understanding of this algebra and its properties it is convenient to consider a
different description of it.

First of all, remind that cobar resolution of a coalgebra $C$ is the tensor
algebra $F(C)=\bigoplus_{i=0}^\infty C^{\otimes n}$, equipped with the
differential
$$\delta_F([c_1|\dots|c_n])=\sum_{i=1}^n(-1)^{i}[c_1|\dots|c_i^{(1)}|c_i^{(2)}|\dots|c_n].$$
Here we have used $|$ instead of the tensor product sign for the
sake of brevity. (Observe that this differential turns $F(C)$ into a
differential graded algebra with respect to the usual free tensor
product in it.)
\begin{prop}
\label{prop2.1} The Weil algebra $W(C)$ of a coalgebra $C$ is
isomorphic to the free differential algebra $\Omega(F(C))$ generated
by algebra $F(C)$. Moreover, the differential $\partial$ of $W(C)$
under this isomorphism is equal to the sum of the free differential
$d$ in $\Omega(F(C))$ and $\delta$ the natural extension of
$\delta_F$ from $F(C)$ to $\Omega(F(C))$.
\end{prop}
\begin{proof}
The isomorphism $\Phi$ in question identifies $i_c$ with $[c]$ and $w_c$ with $d[c]$
(recall that both algebras that we consider are free, so it is enough to define all maps
on their generators). One
easily checks that this map commutes with the differentials. For instance,
$\partial i_c=w_c-i_{c^{(1)}}i_{c^{(2)}}$ is sent by this map to
$d[c]-[c^{(1)}|c^{(2)}]$, which is equal to $(d+\delta)[c]$.
\end{proof}

We shall denote by $W(C)_\natural$ the factor space
$W(C)/[W(C),\,W(C)]$ of $W(C)$ with respect to the commutator
subspace $[W(C),\,W(C)]$ (i.e. the subspace spanned by the {\em
graded\/} commutators
$[\omega_1,\,\omega_2]=\omega_1\omega_2-(-1)^{|\omega_1||\omega_2|}\omega_2\omega_1,\
\omega_1,\omega_2\in W(C)$).

Let now $C$ be a left Hopf-module coalgebra over a Hopf algebra \hc,
i.e. there's an action $\hc\otimes C\to C$ of \hc\ on $C$ satisfying
the relation $(hc)^{(1)}\otimes(hc)^{(2)}=h^{(1)}c^{(1)}\otimes
h^{(2)}c^{(2)}$. Let \mc\ be a left comodule over \hc. Consider the
``crossed-product" $W(C)$\/-bimodule $W(C,\,\mc)=\mc\ltimes W(C)$,
which is isomorphic to $\mc\otimes W(C)$ as vector space, and in
which the right action of $W(C)$ is given by multiplication in right
leg of the tensor product and the left action is defined by the
formula
\begin{equation}
\label{leftact} \alpha\cdot(m\otimes\beta)=m^{(0)}\otimes
S^{-1}(m^{(-1)})(\alpha)\beta.
\end{equation}
In this formula we let \hc\ act on $W(C)$ diagonally, i.e. so that
$h(\alpha\beta)=h^{(1)}(\alpha)h^{(2)}(\beta)$ and we define the action on
the generators as $h(i_c)=i_{h(c)}$ and similarly for $w_c$.

One easily checks that the map $\partial_\mc=1\otimes\partial$ is a differentiation of
the graded bimodule $W(C,\,\mc)$ with respect to the differential $\partial$
on $W(C)$, that is $\partial_\mc(\alpha\cdot\omega)=\partial\alpha\cdot\omega+(-1)^{|\alpha|}\alpha\cdot\partial_\mc\omega$
and $\partial_\mc(\omega\cdot\alpha)=\partial_\mc\omega\cdot\alpha+(-1)^{|\omega|}\omega\cdot\partial\alpha$
for all $\alpha\in W(C)$, $\omega\in W(C,\,\mc)$. Clearly, $\partial_\mc^2=0$, so
$(W(C,\,\mc),\ \partial_\mc)$ and its commutant $W(C,\,\mc)_\natural=W(C,\,\mc)/[W(C,\,\mc),\,W(C)]$
with induced differential $\partial_{\mc,\natural}$ are chain complexes.

Now suppose that \mc\ is a right module and left comodule, verifying
the anti-Yetter-Drinfeld condition (see \eqref{aYDc}). We shall
denote by $W^\hc(C,\,\mc)$ the complex $\mc\otimes_\hc W(C)$ with
differential induced from $W(C)$ as above (we again let \hc\ act on
$W(C)$ diagonally by the same formula as above). In this case the
formula
\begin{equation}
\label{rel1}
m\otimes_\hc\alpha\beta=m^{(0)}\otimes_\hc
S^{-1}(m^{(-1)})(\beta)\alpha
\end{equation}
determines a well-defined equivalence relation on $W^\hc(C,\,\mc)$ and we
denote the resulting factor-complex of $W^\hc(C,\,\mc)$ w.r.t. these relations
by $W^\hc(C,\,\mc)_\natural$ and $\tilde\partial_\mc$ will stand for the
differential in this complex.

Finally consider the ideal $I(C)\subset W(C)$ (the ideal, generated by the curvatures
$w_c$) and its powers. We shall denote by $W_n(C)$ the factor-algebra
$W_n(C)=W(C)/I^{n+1}(C)$ and by $W_n(C)_\natural$ its commuted version
$W_n(C)_\natural=W_n(C)/[W_n(C),W_n(C)]=W_n(C)/[W(C),\,W_n(C)]$. Similarly
$W^\hc_n(C,\,\mc)=\mc\otimes_\hc W(C)/I^{n+1}(C)$ and
$W_n^\hc(C,\,\mc)_\natural$ will denote the factorization of
$W^\hc_n(C,\,\mc)$ by the set of relations, similar to \eqref{rel1}.

\begin{rem}\rm
\label{rem11}
In fact, there's a different way to identify the algebra $W(C)$ with $\Omega(F(C))$.
Namely, consider the map $\Phi'$, sending $i_c$ to $[c]$ and $w_c$ to
$d[c]+[c_{(1)}|c_{(2)}]$ (its inverse sends $[c]$ to $i_c$ and $d[c]$ to $w_c-i_{c_{(1)}}i_{c_{(2)}}$). In this case the differential $\partial$ on $W(C)$
corresponds to $d$ on $\Omega(F(C))$, which is easily checked by a simple
calculation. One direct consequence of this observation is the following
universal property of $W(C)$: {\em for any unital differential algebra $\Omega$ and any
linear map $f$ from $C$ to the space of degree $1$ elements in $\Omega$, there
is a unique map of differential graded algebras $W(C)\to\Omega$, which
coincides with $f$ in degree $1$\/}.
\end{rem}

\subsection{Cohomology of $W^\hc(C,\,\mc)$ and $W^\hc(C,\,\mc)_\natural$}
The identification of proposition \ref{prop2.1} allows one to introduce a
structure of double complex on $W(C)$ simply by transfering it from $\Omega(F(C))$.
This amounts to introducing a bigrading on the free generators $i_c,\ w_c$ of
$W(C)$ so that ${\rm bideg}\,i_c=(0,1)$ and ${\rm bideg}\,w_c=(1,1)$ and presenting
the differential $\partial$ as the sum of two ``partial"
differentials, $\delta$ defined by $\delta(i_c)=-i_{c^{(1)}}i_{c^{(2)}},\
\delta(w_c)=w_{c^{(1)}}i_{c^{(2)}}-i_{c^{(1)}}w_{c^{(2)}}$ and $d,\
di_c=w_c,\ dw_c=0$.

This observation simplifies the calculation of cohomology of $W(C)$ and of its
commutant space $W(C)_\natural$ (compare \cite{Quil}, \S3). Actually,
$W(C)$ being a bicomplex, the
same is true for $W(C)_\natural$ and we can use spectral sequence arguments
to compute the cohomology in both cases. For example, consider the spectral sequence
abutting to the cohomology of  $W(C)$ (resp. $W(C)_\natural$) which starts
from its $d$\/-cohomology. One sees that  $W(C)$ with differential $d$ is
direct sum of tensor powers of the complex $L(C)$, which is equal to $C$ in
dimensions $1$ and $2$ and zero elsewhere with differential equal to the
identity map from $C$ to itself. In the case of  $W(C)_\natural$ tensor powers
are replaced by their ``cyclic" variant, i.e. by tensor powers, factorized
by the action of cyclic groups. As the cohomology of such complex $L(C)$
vanishes, one concludes that $d$\/-cohomologies of  $W(C)$ and  $W(C)_\natural$
vanish, so the spectral sequence collapses and we conclude that the
cohomology of  $W(C)$ and  $W(C)_\natural$ (with respect to the differential
$\partial$) also vanish. Moreover from these considerations one obtains a
contracting homotopy $H$ for the $d$\/-cohomology of $W(C)$, which sends the
elements $w_c$ to $i_c$ and $i_c$ to $0$ and is extended to other elements by the Leibnitz rule.
 Now one can rephrase the
reasoning dealing with $W(C)_\natural$ by saying that $H$ commutes with the
action of cyclic group and hence descends to a contracting homotopy
$H_\natural$ on $W(C)_\natural$.

Similar arguments allow one to calculate the cohomology of the Weil complex of
$C$ in the presence of coefficients' module \mc. It is clear that the bigrading on $W(C)$ induces bigradings on both
complexes $W(C,\,\mc)$ and $W(C,\,\mc)_\natural$ so that the differentials $\partial_\mc$ and $\partial_{\mc,\natural}$
take the form of the sum of two differentials, induce by $\delta$ and $d$.
Now the contracting homotopy $H$ on $W(C)$ above gives rise to homotopies ${}_\mc\!H$ and
${}_\mc\!H_\natural$ of these complexes, so that their $d$\/-cohomologies vanish.
By the spectral sequence argument we conclude that their $\partial_\mc$\/-cohomologies
also vanish. The same argument applies to $W^\hc(C,\,\mc)$ and $W^\hc(C,\,\mc)_\natural$
(to see this it is enough to observe that the contracting homotopy $H$
commutes with the action of \hc, which is clear, since it is completely
determined by its values on the free generators of $W(C)$) so we obtain
the following result:
\begin{prop}
\begin{equation}
H^*(\overline W^\hc(C,\,\mc)_\natural,\,\tilde\partial_\mc)=0,
\end{equation}
where $\overline W^\hc(C,\,\mc)_\natural$ denotes the factor-complex of
$W^\hc(C,\,\mc)$ by the subspace generated by $1\in W(C)$.
\end{prop}

One can also use the isomorphism of remark \ref{rem11} to calculate
the $\partial$\/-cohomology of $W(C)$, $W(C)_\natural$ and
$W^\hc(C,\,\mc)$. Actually this approach makes the arguments
easier, since it allows one to avoid considering the spectral sequence.
In fact, in the case of $W(C)$ and $W^\hc(C,\,\mc)$ this is just the
consequence of the well-known fact that the universal differential
calculus of an algebra is acyclic (presence of coefficients adds
very little to the proof of this statement). And for $W(C)_\natural$
one can use the isomorphism $\Phi'$ and apply the Goodwillie's
(\cite{Good}) and Karoubi's (\cite{Karoubi}) theorems, to show that
$H(W(C)_\natural,\,\partial)\cong H_{DR}(F(C))=0$, where
$H_{DR}(F(C))$ is the non-commutative de Rham cohomology. In other
words, we use the fact that the de Rham cohomology of a graded
algebra without zero degree component vanishes.

However the map $\Phi'$ does not respect the bigrading, hence it is very
hard to understand, where it sends the ideal $I(C)\subset W(C)$ and its
powers. So this identification is not very convenient if one wants to describe the cohomology of
$W_n(C)_\natural$ and $W_n^\hc(C,\,\mc)_\natural$, which is our next
purpose.

But before we can proceed we need to give few definitions, similar to
those given in the paper \cite{HKRS}.

Let $A$ be a left Hopf-module algebra. Dually to the Hopf-type cyclic cohomology of $A$ with
coefficients in \mc\ defined in \cite{HKRS}, we define its Hopf-type cyclic {\em homology\/} as the
homology of the cyclic module $C_\cdot^\hc(A,\,\mc)$, where $C_n^\hc(A,\,\mc)=\mc\otimes_\hc A^{\otimes
n+1}$ and the cyclic operations are induced from the following operations on
$\mc\otimes A^{\otimes n+1}$:
\begin{align}
\partial_i(m\otimes(a_0,a_1,\dots,a_n))&=\begin{cases}
                                        m\otimes(a_0,\dots,a_ia_{i+1},\dots,a_n),\ &0\le i\le n,\\
                                        m^{(0)}\otimes(S^{-1}(m^{(-1)})(a_n)a_1,\dots,a_{n-1}),\ &i=n.
                                        \end{cases}\\
\sigma_i(m\otimes(a_0,a_1,\dots,a_n))&=m\otimes(a_0,\dots,a_{i-1},1,a_i,\dots,a_n)\\
\tau_n(m\otimes(a_0,a_1,\dots,a_n))&=m^{(0)}\otimes(S^{-1}(m^{(-1)})(a_n),a_1,\dots,a_{n-1}).
\end{align}
One easily checks that the anti-Yetter-Drinfeld condition guarantees that
these formulas can be pulled down to $C_\cdot^\hc(A,\,\mc)$ and that if
\mc\ is stable (i.e. $m^{(0)}m^{(-1)}=m$ for all $m\in\mc$) then
$\tau_n^{n+1}=1$. We shall denote the corresponding cyclic homology by
$HC_\ast^\hc(A,\,\mc)$.

Similarly to the coefficient-free case one can define cyclic homology of an
algebra with the help of non-commutative differential formas and various operators on
them. Put $\Omega_\hc(A,\,\mc)=\mc\otimes_\hc\Omega(A)$. As in the
coefficientless case, one defines the Karoubi operator $\kappa$ and the
differentials $b$ and $B$ on $\Omega_\hc(A,\,\mc)$ turning it into a mixed
complex, whose cyclic homology would coincide with $HC_\ast^\hc(A,\,\mc)$.

Also similarly to the coefficient-free case, one can introduce the non-commutative
de Rham homology of $A$ as the homology of the complex $\overline\Omega_\hc(A,\,\mc)$
defined as the factorization of $\Omega_\hc(A,\,\mc)$ modulo the relations similar to
\eqref{rel1}. It is possible to prove a Karoubi-type theorem identifying
this homology with the image of the Connes' $S$\/-operator in $HC_\ast^\hc(A,\,\mc)$.

It is evident that $W(C)_\natural=\overline\Omega(F(C))$ and
$W^\hc(C,\,\mc)_\natural=\overline\Omega_\hc(F(C),\,\mc)$, and one can regard
the reasoning involving the $d$\/-cohomology of $W(C)_\natural$ and
$W^\hc(C,\,\mc)_\natural$ as a particular case of the following simple proposition

\begin{prop}
The non-commutative de Rham cohomology of a free \hc\/-module algebra (including the
case when one introduces coefficients in a anti-Yetter-Drinfeld module) vanishes.
\end{prop}

Observe that in order to define this type of homology and to prove this statement
one does not need the stability conjecture (however if \mc\ is not stable, it is not clear
whether one can interprete this sort of homology \`a la Karoubi's theorem).

For later use we shall need a better knowledge of the structure of the
cyclic homology (with coefficients) of a free algebra.
\begin{prop}
\label{propa}
Let $V$ be a \hc\/-module. Then the Hopf-type Hochschild homology of the
free (tensor) algebra $F=T(V)$ with coefficients in a SAYD module \mc\ vanishes
in dimensions greater than 1. Moreover, the mixed complex associated to
$F$ (see above) is quasi-isomorphic to the following super ($\mathbb Z/2\mathbb Z$\/-graded) complex,
$$
X_\hc(F,\,\mc):\mc\otimes_\hc F\raisebox{1pt}{$\begin{smallmatrix}\stackrel{\scriptstyle{\natural d_\mc}}{\scriptstyle\longrightarrow}\\
\stackrel{\scriptstyle\longleftarrow}{\scriptstyle{
b_\mc}}\end{smallmatrix}$}\Omega_\hc^1(F,\,\mc)_\natural,
$$
and the reduced (i.e. $\mathrm{mod}\,\Bbbk$) cyclic homology of $F$ with arbitrary coefficients is
$$
\overline{HC}^\hc_*(F,\mc)=\begin{cases}
                 (\mc\otimes_\hc\bar F)_\natural,\ &*=0,\\
                  0,\                      &*\ge1
                \end{cases}
$$
\end{prop}
Here $\Omega_\hc^1(F,\,\mc)_\natural=\Omega_\hc^1(F,\,\mc)/\{m\otimes_\hc \omega f-m^{(0)}\otimes_\hc S^{-1}(m^{(-1)})(f)\omega\}$
for all $m\in\mc,\ \omega\in\Omega^1(F),f\in F$ and $(\mc\otimes_\hc\bar F)_\natural$ is defined similarly ($\bar F$ is just $\oplus_{n\ge1}V^{\otimes n}$).
Observe that due to the dimension this space is
equal to $\overline\Omega_\hc^1(F,\,\mc)$. The differentials in this complex are
defined similarly to the coefficientless case: $\natural d_\mc$ is the universal
differential $1_\mc\otimes_\hc d:F\to\Omega_\hc^1(F,\,\mc)$ followed
by the natural projection, and $b_\mc$ is given by the formula
$$
b_\mc(\natural(m\otimes_\hc adb))=m\otimes_\hc ab-m^{(0)}\otimes_\hc
S^{-1}(m^{(-1)})(b)a.
$$
The complex $X_\hc(F,\,\mc)$ is in fact the first level of the Hodge tower
of the mixed complex $(\Omega_\hc(F,\,\mc),\ b,\ B)$, see \cite{CQ3} and section \ref{sct3} below.
\begin{proof}
is essentially the same as that in the absense of coefficients, see e.g.
\cite{Loday}, \S3.1. First one shows that the Hochschild complex of $F$ is
quasi-isomorphic to $C^{\mathrm{small}}(F,\mc)$,
$$
\ldots\longrightarrow0\longrightarrow0\longrightarrow\mc\otimes_\hc(F\otimes
V)\stackrel{b}{\longrightarrow}\mc\otimes_\hc F,
$$
where the nontrivial map $b$ is given by $m\otimes_\hc(f\otimes v)\mapsto m\otimes_\hc fv-m^{(0)}\otimes_\hc
S^{-1}(m^{(-1)})(v)f$ (here and below $v$ will denote an element from $V$ and we omit the tensors, while speaking about $F$).
The quasi-isomorphism is given by the evident inclusion of $C^{\mathrm{small}}(F,\mc)$
into the Hochschild complex of $F$ with coefficients in \mc, $CH^\hc(F,\mc)$ and the projection
$\phi$ from $CH^\hc(F,\mc)$ to $C^{\mathrm{small}}(F,\mc)$ defined
as follows:
$\phi_0=\id$, and
$$
\phi_1(m\otimes_\hc(f\otimes v_1\dots v_n))=
\sum_{i=1}^n m^{(0)}\otimes_\hc (S^{-1}(m^{(-1)})(v_{i+1}\dots
v_n)fv_1\dots v_{i-1}\otimes v_i),
$$
and $\phi_1(m\otimes_\hc(f\otimes 1))=0$.
The contracting homotopy from $CH^\hc(F,\mc)$ to the image of $\phi$ is given by $h_0=0$,
and in higher dimensions we put $h_n(m\otimes_\hc(f_0\otimes\dots\otimes f_{n-1}\otimes
v))=0$ ($v\in V$)and then extend $h_n$ by the recursive formula
\begin{multline*}
h_n(m\otimes_\hc(f_0\otimes\dots\otimes f_{n-1}\otimes f_nv))=\\
h_n(m^{(0)}\otimes_\hc(S^{-1}(m^{(-1)})(v)f_0\otimes\dots\otimes
f_n))+(-1)^nm\otimes_\hc(f_0\otimes\dots\otimes f_n\otimes v).
\end{multline*}
After this one shows that the cyclic bicomplex of $F$ (with coefficients
\mc) is quasiequivalent to
\begin{equation}
\label{com1}
\begin{CD}
       @VVV                          @VVV                     @VbVV\\
         0           @<<< {\mc\otimes_\hc(F\otimes V)} @<\gamma<< {\mc\otimes_\hc F}\\
       @VVV                          @VbVV                    @.   \\
{\mc\otimes_\hc(F\otimes V)} @<\gamma<<{\mc\otimes_\hc F}              @.      {\ }\\
       @VbVV                          @.                      @.   \\
         {\mc\otimes_\hc F}           @.             {\ }               @.      {\ }
\end{CD}
\end{equation}
for the map $\gamma$ given by the equation $\gamma(m\otimes_\hc f)=\phi_1(m\otimes_\hc(1\otimes
f))$. It is an easy exercise to see that this bicomplex is nothing but the
expansion of $X^\hc(F,\,\mc)$. The last statement (concerning the
the reduced homology of a free algebra) follows easily from
the structure of the complex \eqref{com1}
\end{proof}
\begin{rem}\rm
Observe that essentially the same proof can be used to show that the Hochschild
homology (and cohomology) of a free algebra $F=T(V)$ of an \hc\/-module $V$ with
coefficients in
\mc\ {\em and\/} a \hc\/-equivariant bimodule $N$ over $F$
(that is in a bimodule, on which \hc\/ acts so that the actions of \hc\ on $N$ and $F$
agree) vanishes in
degrees exceeding 1. This is a generalization of discussion in \cite{CQ1}
(see the proof that the free algebras are {\em quasi-free\/} in the cited paper).
However we don't know whether it is possible to exclude the explicit
formulas for homotopies etc. from our proof, as it is done in \cite{CQ1}.
This is due to the fact that to our knowledge there's no description of the Hochschild
homology and cohomology in terms of the derived functors in the
\hc\/-equivariant setting (and especially when coefficient module \mc\ appear).
\end{rem}

\subsection{Homology of $W_n^\hc(C,\,\mc)_\natural$ and $S$\/-operations}
Let us now denote the cohomology of the complex $W_n^\hc(C,\,\mc)_\natural$
by $H^*_\hc(C,\,\mc;\,n)$. Note that obviously $H^*_\hc(C,\,\mc;\,0)\cong
HC^{*-1}_\hc(C,\,\mc)$. Now the following statement is a straightforward generalization of the Theorem 7.1 in
\cite{Crain}.
\begin{theorem}
\label{theoa}
There are canonical isomorphisms $\alpha_n:H^*_\hc(C,\,\mc;\,n)\cong
HC^{*-1-2n}_\hc(C,\,\mc),\ n\ge0$. Moreover, under this identification the homomorphism
 induced by the canonical projection $\pi_n:W_n^\hc(C,\,\mc)_\natural\to
W_{n-1}^\hc(C,\,\mc)_\natural$ coincides with the $S$\,-operation in the
cyclic cohomology.
\end{theorem}
\begin{proof}
The proof can be obtained as a word to word repetition of that in
\cite{Crain}. Here we shall give a slightly modified version of this
reasoning.

First of all, we observe that similarly to the cited paper, in our case it
is enough to prove the non-\hc\/-equivariant, coefficientless version of this
theorem. To see this, just notice that all the statements below commute
with the \hc\/-action. (However for the convenience of the reader in the end we shall
show, what is to be changed in the general case.)

Further we consider the case $n=0$. As it is observed in the last sentence preceeding
this theorem, the statement holds in this particular case (actually $W_0(C)_\natural$ is isomorphic to
$C_\ast^\lambda(C)[1]$, and hence $H^*(C;\,0)=HC^{*+1}(C)$).

In order to proceed from $n=0$ to $n\ge1$, we use the following proposition, similar to the Lemma
8.2 of \cite{Crain}:
\begin{lemma}
\label{lema}
Let $I^{(n)}_\natural(C)$ denote the subspace of $W(C)_\natural$ spanned by
the elements containing exactly $n$ factors of type $w_{c}$. Then
\begin{description}
\item[{\rm({\em i\/})}] there is isomorphism
$p:H^*(W_n(C)_\natural)\stackrel{\cong}{\to}H^*(I^{(n)}_\natural(C)/Im\/d,\,\delta)$,
compatible with $S$\,-operation;
\item[{\rm({\em ii\/})}] the isomorphism $p$ identifies the map $\pi^*_n:H^*(W_n(C)_\natural)\to
H^*(W_{n-1}(C)_\natural)$ with the coboundary operation in the long exact sequence,
associated to the short exact sequence
\begin{equation}
\label{eqb}
0\longrightarrow(I^{(n-1)}_\natural(C)/Im\/d,\delta)
\stackrel{d}{\longrightarrow}
(I^{(n)}_\natural(C),\delta)\longrightarrow(I^{(n)}_\natural(C)/Im\/d,\delta)\longrightarrow0,
\end{equation}
(observe that the dimension of an element is shifted by $d$)
i.e. the following diagram commutes
$$
\begin{CD}
{H^*(W_n(C)_\natural)}                   @>\pi_n^*>>              {H^*(W_{n-1}(C)_\natural)}\\
 @VpVV                                                               @VpVV\\
{H^*(I^{(n)}_\natural(C)/Im\/d,\,\delta)}@>{\delta_n^*}>>{H^*(I^{(n-1)}_\natural(C)/Im\/d,\,\delta)}.
\end{CD}
$$
\end{description}
\end{lemma}
This lemma reduces our theorem to an equivalent statement about the complex
$(I^{(n)}_\natural(C)/Im\/d,\delta)$. Now we recall that $W(C)$ is
isomporphic to $\Omega(F(C))$ --- the universal differential calculus of the
cobar-resolution of $C$, and that the differential $d$ in $W(C)$ corresponds
to the universal differential in $\Omega(F(C))$. So, we are to compare the
cyclic coalgebra
cohomology of $C$, $HC^*(C)$ with the homology of the complex $(\overline{\Omega}^n(F(C))/Im\/d,\,\delta)$.

To this end we recall from \S 2.6 \cite{Loday} and \cite{Karoubi}, \S 2.12-2.14, that the space
$\overline{\Omega}^n(A)/Im\/d$, which is the same as $\overline C_n(A)_{\mathrm{ab}}/Im\/d$
 in the cited paper (here $\overline C_n(A)$ denotes the normalized cyclic complex
 of an algebra $A$, i.e. complex, calculating the $\mathrm{mod}\,\Bbbk$ homology of $A$)
 is isomorphic to $\overline
C_n^\lambda(A)/Im\/b$ ($b$ is the Hochschild differential). Thus we are to
compare the homology of $(\overline C_n^\lambda(F(C))/Im\/b,\,\delta)$, where $\delta$ is the intrinsic differential in $F(C)$, with
$HC^*(C)$.

So suppose $n\ge1$ and consider the complex $(\overline
C_n^\lambda(F(C))/Im\/b,\,\delta)$. We shall define an isomorphism $\varphi:H^*(\overline
C_n^\lambda(F(C))/Im\/b,\,\delta)\stackrel{\cong}{\to}H^{*+2}(\overline
C_{n+1}^\lambda(F(C))/Im\/b,\,\delta)$.
Let $[x]$ be an element in its homology,
represented by a (co)cycle $x\in\overline C_n^\lambda(F(C))/Im\/b$. We can
choose an element $x'$ in $\overline C_n^\lambda(F(C))$, equal to $x$ modulo
the image of $b$. Then since $x$ is a cocycle in $(\overline
C_n^\lambda(F(C))/Im\/b,\,\delta)$, we conclude that $\delta x'\in Im\/b$. Choose $w'\in\overline
C_{n+1}^\lambda(F(C))$ so that $\delta x'=bw'$. Let $w$ be the image of $w'$
under the natural projection $\overline C_{n+1}^\lambda(F(C))\to\overline
C_{n+1}^\lambda(F(C))/Im\/b$. Then (all the equalities are strict and not
$\mathrm{mod}\,b$ since $b^2=0$)
\begin{equation*}
b(\delta w)= b(\delta w')=\delta(bw')=\delta\delta x'=0.
\end{equation*}
Thus $\delta w\in Ker\/b$. But the complex $\overline
C_*^\lambda(F(C))$ is acyclic in dimensions $*\ge1$ (see proposition \ref{propa}), so
$Ker\/b=Im\/b$, and we see that $\delta w=0$. Put $\varphi([x])=[w]$.

Let us show, that $\varphi$ is well-defined. There were three ambiguities
in its definition. First when we chose a representative cocycle $x$
in $[x]$, second when we passed from $x$ to $x'$ and finally, when we
found $w'$. Now if $\bar x$ is a different representative of $[x]$ and $\bar x'$ is
its preimage in $\overline C_n^\lambda(F(C))$, then $x'-\bar x'=\delta\alpha
+ b\beta$ for some $\alpha\in\overline C_n^\lambda(F(C)),\ \beta\in\overline
C_{n+1}^\lambda(F(C))$. Then
\begin{equation}
\label{eqa}
\delta x'-\delta\bar
x'=\delta(\delta\alpha+b\beta)=\delta b\beta=b\delta\beta.
\end{equation}
Now let $\bar w'$ verify the equation the $b\bar w'=\delta\bar x'$. It
follows from \eqref{eqa}, that $b(w'-\bar w')=b\delta\beta$. From the
acyclicity of $\overline C_*^\lambda(F(C))$ with respect to $b$ it
follows that $w'-\bar w'=\delta\beta+b\gamma$ for some $\gamma\in\overline
C_{n+2}^\lambda(F(C))$. Hence the difference $w-\bar w$, where $\bar w$ is
the image of $\bar w'$ under the natural projection ${\rm mod}\,Im\/b$, is equal to the image
of $\delta\beta$ under this projection. So $[w]=[\bar w]$.

In order to show, that $\varphi$ is an isomorphism, we shall construct its
inverse $\psi$. To this end observe that since the differential graded (nonunital) algebra
$F(C)$ is acyclic, so are all the complexes $(\overline C_n(F(C)),\,\delta)$
and $(\overline C_n^\lambda(F(C)),\,\delta)$ for all $n\in\mathbb N$. Now we can define
the map $\psi$ as follows. Take a representative $y$ of $[y]\in H^*(\overline
C_{n+1}^\lambda(F(C))/Im b,\,\delta)$, then $by\in\overline
C_n^\lambda(F(C))$ is a well-defined element (since $b^2=0$) and $\delta by=-b\delta
y=0$, since $y$ is a cocycle. By the virtue of acyclicity of $(\overline
C_n^\lambda(F(C)),\,\delta)$ there exists an element $v'\in\overline
C_n^\lambda(F(C))$, such that $\delta v'=by$. Thus $v'$ is a $\mathrm{mod}Im\/b$ cocycle in $(\overline
C_n^\lambda(F(C)),\,\delta)$ and we put $\psi([y])=[v]$, where $v\in\overline C_n^\lambda(F(C))/Im\/b$ is the projection of $v'$. The proof
that $\psi$ is well-defined is as easy as before.

Now we can show that $\psi\circ\varphi=\varphi\circ\psi=1$. In fact if
$[w]=\varphi([x])$, then $bw=bw'=\delta x'$ and by the very definition of
$\psi$ $\psi([w])=[x]$ (since $x'$ projects in $x$ modulo $Im\/b$).
Vice-versa, if $\psi([y])=[v]$, then $\delta v'=by$ (where $v'$ is the
preimage of the cocycle $v$ representing $[v]$ under the natural projection)
and, again by the definition of $\varphi$, $\varphi([v])=[y]$. Thus
$\varphi$ is an isomorphism and we put $\alpha_n=p^{-1}\circ\varphi^n$.

Now it only remains to show, that the iterations of $\varphi$ send
$S$\/-operations on the cyclic cohomology of a coalgebra to the projection
$\pi_n$, or, more accurately, that
$$
S=\sigma^{-1}\circ\psi^{n-1}\circ p\circ\pi^*_n\circ
p^{-1}\circ\varphi^n\circ\sigma.
\footnote{This is not quite true: below we shall have to introduce some scalar
coefficients into the game in order to make the diagram commute.}
$$
Here $\sigma:H^*(C)\stackrel{\cong}{\to} H^{*+1}(W_0(C)_\natural)$ is
the natural isomorphism induced by the dimension shift and $p$ is the isomorphism of the lemma \ref{lema}. In the future we
shall omit the ``suspension" $\sigma$, since it doesn't change the proof, except for,
probably, the introduction of some signs.

Once again we start with the case $n=1$. We have to compare the
morphism $S$ with $p\circ\pi_1\circ p^{-1}\circ\varphi$. First of all we
recall from ({\em ii\/}) Lemma \ref{lema} that $p\circ\pi_1\circ
p^{-1}=\delta^*_1$, the coboundary operation of the exact sequence
\eqref{eqb} for $n=1$, so we are reduced to showing that
$S=\delta^*_1\circ\varphi$.

Consider the exact sequence, which is used to define the $S$\/-operation in
cohomology of $F(C)_\natural$:
\begin{equation}
\label{eqc}
0\longrightarrow F(C)_\natural\stackrel{d}{\longrightarrow}\Omega^1(F(C))_\natural\stackrel{b}{\longrightarrow}F(C)\stackrel{\natural}{\longrightarrow}F(C)_\natural\longrightarrow
0,
\end{equation}
where, for an algebra $R$,
$\Omega^1(R)_\natural\eqdef\Omega^1(R)/[\Omega^1(R),\,R]$ (this space was denoted $\overline\Omega^1(R)$ above).
This sequence is exact, because $F(C)$ is a free algebra, and the map $S$ is
defined as the result of the diagram chasing: for a cocycle $x\in F(C)_\natural$ one
finds $x'\in F(C)$ that surjects on $x$, then $\natural\delta(x')=0$ and hence
$\delta x'=by$. Once again $b\delta y=0$ and one can find $z$ in $F(C)_\natural$ such that
$dz=\delta y$. The map $S$ is then given by $x\mapsto z$.
Observe now that the sequence \eqref{eqc} is the result of splicing
together two short exact sequences:
\begin{align}
\label{eqd}
&0\longrightarrow\Omega^1(F(C))_\natural/Im\/d\longrightarrow
F(C)\stackrel{\natural}{\longrightarrow}F(C)_\natural\longrightarrow0\\
\intertext{and}
\label{eqe}
0&\longrightarrow
F(C)_\natural\stackrel{d}{\longrightarrow}\Omega^1(F(C))_\natural\longrightarrow\Omega^1(F(C))_\natural/Im\/d\longrightarrow0.
\end{align}
So the map $S$ is equal to the composition of the coboundary operations in
the long exact sequences associated to \eqref{eqd} and \eqref{eqe}.

But $I^{(0)}\cong F(C)$,
$I^{(1)}=\Omega^1(F(C))$, and consequently  $I^{(0)}_\natural=F(C)_\natural$,
$I^{(1)}_\natural=\Omega^1(F(C))_\natural$, $I^{(1)}_\natural/Im\/d=\Omega^1(F(C))_\natural/Im\/d=\overline C_1^\lambda(F(C))/Im\/b$. Now it is evident,
that $\varphi$ is equal to the coboundary operation of \eqref{eqd} (remark that it is
obviously an isomorphism because $F(C)$ is acyclic) and
$\delta_1^*$ --- to that of \eqref{eqe}. The proof in case $n=1$ is finished.

The general case is reduced to the one we have just considered by the
following observations. First, the Hochschild homology of $F(C)$ vanishes in dimensions $\ge2$ and hence the maps $\varphi$ and $\psi$ can be extended from
$I^{(n)}_\natural/Im\/d=\overline C^\lambda_n(F(C))/Im\/b$ to $I^{(n)}_\natural\cong\Bigl(\overline
C_n(F(C))/Im\/b\Bigr)/Im(1-t)$ where $t$ is the cyclic operator (see
\cite{Karoubi}, \S 2.12-2.16) when $n\ge1$.

Now an easy diagramm chasing shows
that (up to a scalar multiple) these extensions will commute with the maps in the exact sequences
\eqref{eqb} for $n$ and $n+1$ respectively. To this end consider the following  diagram
\begin{equation*}
\begin{CD}
0 @>>> I^{(n)}_\natural/Im d @>d>> I^{(n)}_\natural @>b>> \widehat{ \Omega}^n(F(C))@>>>
I^{(n)}_\natural/Im d @>>> 0\\
@.       @VV\displaystyle{\cdot \frac n{n+1}}V  @VV\displaystyle{-sb}V       @VV\displaystyle{sd}V     @|       @.\\
0 @>>> I^{(n)}_\natural/Im d @>b>>\widehat{ \Omega}^{n-1}(F(C)) @>d>>
I^{(n)}_\natural @>>> I^{(n)}_\natural/Im d @>>>
0.
\end{CD}
\end{equation*}
Here we identify $I^{(n)}_\natural/Im\/d$ with $\Omega^n(F(C))/Im\/d$
 and $I^{(n)}_\natural$ with $\Omega^n/Im\,b +Im(1-\kappa)$
(see \cite[\S 3]{CQ2} for reference), $\widehat{ \Omega}^{n}(F(C))$ is an abbreviation for $ \Omega^n(F(C))/Imd+Im(1-\kappa)$ and $s$ is the standard
contracting homotopy for the universal differential $d$ of $\Omega(F(C))$.
This diagram is commutative because of the identities
\begin{gather*}
-sbd = -\frac 1{n+1}sbdN_n = -\frac 1{n+1}sbB = \frac 1{n+1} sBb =
\frac 1{n+1} sdN_{n-1}b =\\
\shoveright{\frac n{n+1}(b-dsb) = \frac n{n+1}b,}\\
\shoveleft{sdb = (\id-ds)b = b-dsb  = d(-sb),}\\
\shoveleft{sd = \id - ds = \id.}
\end{gather*}

It follows that we have a commutative diagram for boundary homomorphisms. The
boundary morphism of the top row is $\pi\circ\varphi$ and that of the bottom row
is equal to $\varphi\circ\pi$. So we obtain the equality
$$
\frac n{n+1}\pi_{n+1}\circ\varphi_n = \varphi_{n-1}\circ\pi_n.
$$
Thus we see that $\varphi$ commutes with $\pi$. Since $\psi$ is the inverse of
$\varphi$, the same is true for $\psi$.

 So we conclude that $\varphi$, rescaled, if necessary, so as to eliminate the factor  $\dfrac{n}{n+1}$ is
an isomorphism of the exact sequences ($\varphi^{-1}=\psi$, or its rescaling), in particular,
$\delta^*_{n+1}=\varphi\circ\delta^*_n\circ\psi$. So,
\begin{equation*}
\begin{split}
\sigma^{-1}\circ\psi^{n-1}\circ p\circ\pi^*_n\circ
p^{-1}\circ\varphi^n\circ\sigma&=\sigma^{-1}\circ\psi^{n-1}\circ\delta^*_n\circ\varphi^n\circ\sigma\\
                               &=\sigma^{-1}\circ\psi^{n-2}\circ\delta^*_{n-1}\circ\varphi^{n-1}\circ\sigma\\
                               &=\dots=\\
                               &=\sigma^{-1}\circ\delta^*_1\circ\varphi\circ\sigma=S
\end{split}
\end{equation*}

Finally, let us briefly describe, what changes should be made in this
reasoning to make it work in the general situation, when the coefficients
appear.

First of all, observe that an analog of lemma \ref{lema} holds with $I^{(n)}_\natural(C)$
replaced with $I^{(n)}_\hc(C,\,\mc)_\natural$ the subspace of $W_\hc(C,\,\mc)_\natural$ spanned by the elements, containing exactly
$n$ factors of type $w_{c}$. Further, the isomorphism $\overline C_n(A)_{\mathrm{ab}}/Im\/d=\overline
C_n^\lambda(A)/Im\/b$ can be extended to $\overline C^\hc_n(A,\,\mc)_{{\mathrm
ab}}/Im\/d=(\overline C_n^\hc)^\lambda(A,\,\mc)/Im\/b$ for all \hc\/-module algebras $A$.
Finally, observe that the properties of the cyclic complexes used in definitions of $\varphi$ and $\psi$ and in the proof of the $S$\/-operation relation, such as the acyclicity of the complexes $(\overline
C_*^\lambda(F(C)),\,b)$ and $(\overline C_n^\lambda(F(C)),\,\delta)$, hold
for $((\overline C^\hc_*)^\lambda(F(C),\,\mc),\,b)$ and $((\overline
C^\hc_n)^\lambda(F(C),\,\mc),\,\delta)$.
\end{proof}

\begin{rem}\rm
In general, the following diagram shows that the composition $\pi\circ\varphi$
coincides (up to scalar multiplier) with $S$-operation described in
~\cite{Crain}
\begin{equation*}
\begin{CD}
0 @>>> \Omega^n/b,d @>d>> \Omega^{n+1}/b,(1-\kappa) @>b>>  \Omega^n/d,(1-\kappa) @>>> \Omega^n/b,d @>>> 0\\
@.       @|  @VV\displaystyle{sN}V       @VV\displaystyle{N}V     @VV\displaystyle{\cdot n}V       @.\\
0 @>>> \Omega^n/b,t-1 @>\mathcal{N}>> \Omega^{n}/d @>t-1>>
\Omega^n/d @>>> \Omega^n/d,t-1 @>>> 0.
\end{CD}
\end{equation*}
Here $t$ is a cyclic permutation on $W(C)$
\begin{equation*}
t(a\omega) = (-1)^{|a||\omega|}\omega a,\quad \omega\in W(C),\ a=i_c\mbox{ or } a=w_c,\
c\in C,
\end{equation*}
and $\nc(x) = \sum_{i=0}^p t^i(x)$, where $p$ is $\delta$-degree of
the element $x$.
\end{rem}

\begin{rem}\rm
\label{rem1}
In future we shall need to know not only the homology of
$W^\hc_n(C,\,\mc)_\natural$, but also that of $I^\hc_n(C,\,\mc)_\natural$,
which is equal to $I^\hc_{n+1}(C,\,\mc)/[I^\hc_{n}(C,\,\mc),\,I(C)]$.Here
$I^\hc_{n+1}(C,\,\mc)$ (resp. $I^\hc_{n}(C,\,\mc)$) is
the subspace of $W^\hc(C,\,\mc)$, generated by the $n+1$\/-st (resp. $n$\/-th) power
of the ideal $I(C)$. But similarly to the proof of theorem 6.7 of \cite{Crain}
one shows first that $I^\hc_n(C,\,\mc)_\natural$ is quasiisomorphic to
$\tilde I^\hc_n(C,\,\mc)_\natural=I^\hc_{n+1}(C,\,\mc)/I^\hc_{n+1}(C,\,\mc)\cap[W^\hc(C,\,\mc),\,W(C)]$
and from the long exact sequence of the three-term exact
sequence
$$
0\longrightarrow\tilde I^\hc_{n+1}(C,\,\mc)_\natural\longrightarrow
W^\hc(C,\,\mc)_\natural \longrightarrow
W^\hc_n(C,\,\mc)_\natural\longrightarrow 0
$$
where in the middle stands an acyclic complex, we find that $H^*(I^\hc_{n+1}(C,\,\mc)_\natural)\cong
H^*(\tilde I^\hc_{n+1}(C,\,\mc)_\natural)\cong H^{*-1}(W^\hc_n(C,\,\mc)_\natural)=H^{*+1}_\hc(C,\,\mc;\,n)\cong
HC^{*-2-2n}_\hc(C,\,\mc)$. We shall denote the resulting isomorphism by
$\beta_n$.

Observe, that although the
maps $\alpha_n,\ \beta_n$ and their inverses are defined only on the level
of cohomology, one can write down explicit formulas for them on the level of
the chain complexes. To this end one shall need to choose the contracting
homotopies of all the acyclic complexes that appeared in the proof and use an
explicit isomorphism $\overline{\Omega}^n(A)/Im\/d\cong\overline
C_n^\lambda(A)/Im\/b$. This is precisely what is done in the the section 8
of \cite{Crain} with the help of the maps $\varphi_0,\ \varphi_1,\theta$
etc..
\end{rem}

There is one more way to define the $S$\/-operations in terms of the Weil
complex cohomology. To this end consider the $X$\/-complex of $W(C)$ or
rather the $X_\hc$\/-complex with coefficients \mc, see above. Operations on this
super-complex commute with the differential $\partial$ on $W(C,\,\mc)$, so
by the same computations as in Corollaries 6.9 and 6.10 and theorem 7.9 of \cite{Crain}
we obtain the following statements:
\begin{prop}
\label{longcom}
For all $n$ there are long exact sequences of complexes:
\begin{align}
CC&(W_n^\hc(C,\,\mc)):\dots\lar W^\hc_n(C,\,\mc;b)\stackrel{t-1}{\lar}W_n^\hc(C,\,\mc)\stackrel{N}{\lar}W^\hc_n(C,\,\mc;b)\lar\dots\\
\label{comw}
0&\lar W^\hc_n(C,\,\mc)_\natural\stackrel{N}{\lar}W^\hc_n(C,\,\mc;b)\stackrel{t-1}{\lar}W_n^\hc(C,\,\mc)\stackrel{N}{\lar}W^\hc_n(C,\,\mc;b)\lar\dots\\
\label{comi}
0&\lar\,\,
I^\hc_n(C,\,\mc)_\natural\,\stackrel{N}{\lar}\,\,I^\hc_n(C,\,\mc;b)\,\stackrel{t-1}{\lar}\,\,I_n^\hc(C,\,\mc)\,\stackrel{N}{\lar}\,\,I^\hc_n(C,\,\mc;b)\,\lar\,\,\dots.
\end{align}
Here we use the notation of previous remark and the following agreements: $t$ is the twisted cyclic permutation on $W^\hc(C,\,\mc)$, $t(m\otimes xa)=(-1)^{|a||x|}m^{(0)}\otimes S^{-1}(m^{(-1)})(a)x$ for
all $x\in W(C)$ and generator $a$ of $W(C)$; $N=1+t+t^2+\ldots+t^{p-1}$ for
all elements of tensor degree $p$; and $W^\hc_n(C,\,\mc;b)$ is the complex
$W^\hc_n(C,\,\mc)$ with differential $b=\partial+b_t$,
$$
b_t(m\otimes xa)=(-1)^{|a|}t(a\delta(x)).
$$
In particular the bicomplexes \eqref{comw} and \eqref{comi} compute the
cohomology of $W^\hc_n(C,\,\mc)_\natural$ and $I^\hc_n(C,\,\mc)_\natural$
respectively. The $S$\/-operation in the cohomology of $W^\hc_n(C,\,\mc)_\natural$ and $I^\hc_n(C,\,\mc)_\natural$
is given by the shift of these complexes.
\end{prop}
\begin{prop}
\label{shortcom}
In the previous notation, there exist short exact sequences of complexes
\begin{align}
\label{comw1}
0&\lar W^\hc_n(C,\,\mc)_\natural\stackrel{N}{\lar}W^\hc_n(C,\,\mc;b)\stackrel{t-1}{\lar}W_n^\hc(C,\,\mc){\lar}W^\hc_n(C,\,\mc)_\natural\lar0\\
\label{comi1}
0&\lar\,\,
I^\hc_n(C,\,\mc)_\natural\,\stackrel{N}{\lar}\,\,I^\hc_n(C,\,\mc;b)\,\stackrel{t-1}{\lar}\,\,I_n^\hc(C,\,\mc)\,{\lar}\,\,I^\hc_n(C,\,\mc)_\natural\,\lar\,\,0.
\end{align}
The $S$\/-operation in the cohomology of $W^\hc_n(C,\,\mc)_\natural$ and
$I^\hc_n(C,\,\mc)_\natural$ is given by the cup-product with the
$Ext^2$\/-class, determined by the sequences \eqref{comw1}, \eqref{comi1}.
\end{prop}
It is easy to see that both methods give the same result, i.e. the
$S$\/-operation doesn't depend on the way we define it. To see this,
observe, that the sequence \eqref{comw1} is the version of
\eqref{eqc} with $F(C)$ replaced with $W_n(C)$ (for a while we omit the coefficients
from our notations). As in the lemma \ref{lema} we can replace $W_n(C)$ by
$I^{(n)}(C)$ in this exact sequence. Then as we have shown in the end
of the proof of theorem \ref{theoa}, the cup product with the class of this
exact sequence is equal to the composition of the map $\varphi$ and
$\delta^*_n$, which we have identified with $S$.

\section{Pairings}
The purpose of this section is to define the pairing of $HC^*_\hc(C,\,\mc)$
and $HC^*_\hc(A,\,\mc)$ ($A$ is an \hc\/-module algebra on which $C$ acts on the left in a way,
described in the section 1). Our approach to this
question will be again a suitable generalization of the one used in paper
\cite{Crain}. In the next section we shall compare this
construction with a construction, generalizing the traditional way of
introducing multiplications and comultiplications in cyclic homology.

\subsection{Equivariant $X$\/-complexes with coefficients}
\label{sct3}
Before we define the pairing of the cyclic cohomologies of algebras and
coalgebras we shall briefly discuss the way one generalizes the
Cuntz-Quillen tower of $X$\/-complexes to embrace the cyclic homology and
cohomology with coefficients. Below we shall need the description of cyclic
cohomology of an \hc\/-module algebra based on these ideas.

First of all we shall work in the category of (left) \hc\/-module algebras
$\ac^\hc$ and their morphisms. If $A\in Ob\ac^\hc$ and \mc\ is a (stable)
anti-Yetter-Drinfeld module over \hc, then one can
define the {\em universal differential calculus of $A$ with values in \mc\/}
as it is done in previous section and the Hochschild and cyclic operators on it. As before,
we shall denote
the module of $\mc$-valued differential forms by $\Omega_\hc(A,\,\mc)$ and consider
the operators
$b$ and $B$. One defines the {\em Hodge filtration\/} on $\Omega_\hc(A,\,\mc)$ and
the associated {\em Hodge tower\/} of supercomplexes
by the same formulas as in \cite{CQ2}:
\begin{align}
F^n\Omega_\hc(A,\,\mc)&=b\Omega^{n+1}_\hc(A,\,\mc)\oplus\bigoplus_{k>n}
\Omega^k_\hc(A,\,\mc)\\
\theta\Omega_\hc(A,\,\mc)&=(\Omega_\hc(A,\,\mc)/F^n\Omega_\hc(A,\,\mc)).
\end{align}

The tower $\theta\Omega_\hc(A,\,\mc)$ is a {\em special tower\/}, and
similar to the coefficientless
case we have the following identities
 (here $\Omega=\Omega_\hc(A,\,\mc)$)
\begin{equation}
\label{eqh}
\begin{split}
H_\nu(F^{n-1}\Omega/F^n\Omega)&=\begin{cases}
                                       HH^\hc_n(A,\,\mc), &\nu=n+2\mathbb Z,\\
                                       0,                 &\nu=n-1+2\mathbb Z,
                                \end{cases}\\
H_\nu(\Omega/F^n\Omega)&=       \begin{cases}
                                       HC^\hc_n(A,\,\mc), &\nu=n+2\mathbb Z,\\
                                       HD^\hc_n(A,\,\mc), &\nu=n-1+2\mathbb Z,
                                \end{cases}\\
H_\nu(\widehat\Omega)&=        HP_\nu^\hc(A,\,\mc).
\end{split}
\end{equation}
Here $HH_n,\ HC_n,\ HD_n$ and $HP_\nu$ are respectively Hochschild, cyclic,
(non commutative) de Rham homology and the $\mathbb Z/2\mathbb Z$ graded periodic
homology
of $A$, and $\widehat\Omega\eqdef{\displaystyle{\lim_{\longleftarrow}}}\,\Omega/F^n\Omega$.
Similarly, one retrieves the cyclic (periodic, Hochschild) {\em cohomology\/} of $A$ (with
coefficients in \mc) by considering the dual complex
$\Hom(\theta\Omega,\,k)$. The first level $\Omega/F^1\Omega$ of the Hodge tower
is called the $X$\/-complex of $A$ (with
coefficients in \mc). It is precisely the super-complex $X_\hc(A,\mc)$
that appeared in the proposition \ref{propa}.

Similarly to the cited paper of Cuntz and Quillen one can define a filtration on
$X_\hc(R,\,\mc)$ associated to an ideal $I$ in algebra $R$ ($I$
should be stable under the action of
\hc):
\begin{equation}
\label{filtr}
\begin{split}
&F^{2n+1}_IX_\hc(R,\,\mc):\mc\otimes_\hc
I^{n+1}\leftrightarrows\natural(\mc\otimes_\hc(I^{n+1}dR+I^ndI)),\\
&F^{2n}_IX_\hc(R,\,\mc):\mc\otimes_\hc I^{n+1}+[\mc\otimes_\hc I^n,\,R]\leftrightarrows\natural(\mc\otimes_\hc(I^ndR)).
\end{split}
\end{equation}
In this formulas and below we denote by $[\mc\otimes_\hc I^n,\,R]$ and similar
expressions the image of the map $[,]:\mc\otimes_\hc(I^n\otimes R)\to\mc\otimes_\hc R$
given by the formula
\begin{equation}
\label{commute}
[,](m\otimes_\hc(x\otimes r))=m\otimes_\hc(xr)-m^{(0)}\otimes_\hc(S^{-1}(m^{(-1)})(r)x).
\end{equation}
for all $m\in\mc,\ x\in I^n,\ r\in R$. This map is well-defined, because \mc\ is a
AYD-module. Similarly the subscript $\natural$ or the same sign in front of an
expression denotes the factorization by the space of commutators.

Then we put $\xc_\hc^p(R,\,I;\,\mc)=X_\hc(R,\,\mc)/F^p_IX_\hc(R,\,\mc)$:
\begin{equation}
\begin{split}
&\xc_\hc^{2n+1}(R,\,I;\,\mc):\mc\otimes_\hc(R/I^{n+1})\leftrightarrows\Omega^1_\hc(R,\,\mc)_\natural/\natural(\mc\otimes_\hc(I^{n+1}dR+I^ndI)),\\
&\xc_\hc^{2n}(R,\,I;\,\mc):\mc\otimes_\hc R/\mc\otimes_\hc I^{n+1}+[\mc\otimes_\hc
I^n,\,R]\leftrightarrows\Omega^1_\hc(R,\,\mc)_\natural/\natural(\mc\otimes_\hc(I^ndR)).
\end{split}
\end{equation}
We apply
this construction to $R=RA\eqdef\Omega^{even}(A)={\displaystyle\bigoplus_{i\ge0}}\Omega^{2i}(A)$,
 equipped with {\em Fedosov product\/},
$\omega\circ\omega'=\omega\omega'+(-1)^{|\omega|}d\omega d\omega'$, and the ideal
$I=IA={\displaystyle\bigoplus_{i>0}}\Omega^{2i}(A)$ (the action of \hc\ on
$RA$ is defined via its action on $A$). Thus we obtain a tower
$\mathcal X_\hc(A,\,\mc)=\mathcal X_\hc(RA,\,IA;\,\mc)=(\mathcal
X_\hc^n(RA,\,IA;\,\mc))$ of supercomplexes, which is again a special tower
and its homology (resp. cohomology) verifies the same equations \eqref{eqh}. Moreover this
tower is homotopy equivalent to $\theta\Omega_\hc(A,\,\mc)$. So the
cyclic-type homology of $A$ (with coefficients in \mc) is given by the formulas,
similar to \eqref{eqh} with $\xc_\hc(A,\,\mc)$ instead of $\theta\Omega$. In
particular
\begin{align*}
H_\nu(\xc_\hc^p(A,\,\mc))&=       \begin{cases}
                                       HC^\hc_n(A,\,\mc), &\nu=n+2\mathbb Z,\\
                                       HD^\hc_n(A,\,\mc), &\nu=n-1+2\mathbb Z,
                                \end{cases}\\
H_\nu(\widehat{\xc_\hc}(A,\,\mc))&=        HP_\nu^\hc(A,\,\mc),
\end{align*}
where
$\widehat{\xc_\hc}(A,\,\mc)={\displaystyle{\lim_{\longleftarrow}}}\,\xc_\hc^p(A,\,\mc)$.
Further one can extend the definition of {\em quasi-free\/} algebras given in \cite{CQ1,CQ2}
to the category
$\ac^\hc$, by saying that algebras in this category are quasi-free, if they
verify the conditions,
similar to those, listed in \cite{CQ1} Prop. 3.3 (see also \cite{CQ2} Prop. 7.1),
only this time all the morphisms should be in $\ac^\hc$ (the only problem is
with the condition, concerning the homology dimension with respect to the
Hochschild cohomology, see remark following the proposition \ref{propa}).
Then one can show, that for any exact sequence of \hc\/-algebras with quasi-free $\rc$
\begin{equation}
\label{eqseq}
0\longrightarrow I\longrightarrow R\longrightarrow A\longrightarrow0,
\end{equation}
which is splittable as a sequence of \hc\/-modules (that is we suppose that
one can choose a linear splitting $\rho:A\to\rc$, such that
$\rho(h(a))=h(\rho(a))$), the tower $\xc_\hc(R,\,I;\,\mc)$ is equivalent
to $\xc_\hc(A,\,\mc)$. So one can use {\em any\/} quasi-free extension of $A$ to
calculate cyclic (co)homology of $A$. The extension
$$ 0\longrightarrow IA\longrightarrow RA\longrightarrow A\longrightarrow0$$
is the
{\em universal quasi-free extension\/} of $A$
(this means that for any extension \eqref{eqseq} of $A$ there is a map
$(RA,\,IA)\to(R,\,I)$ covering the identity map on $A$).

\subsection{Hopf-type cyclic cohomology and higher traces}
In view of the statements formulated in the section \ref{sct3}, one can give the following description of the
cyclic cohomology of an algebra $A$. First of all, as we have already said,
$H^\nu(\xc_\hc^p(R,\,I;\,\mc))=HC_\hc^p(A,\,\mc),\ \nu=p+2\mathbb Z$ for any
quasi-free extension \eqref{eqseq} of $A$. Thus for $p=2n$ all the
cohomology classes of $A$ are representible by linear functionals $\mc\otimes_\hc R/\mc\otimes_\hc I^{n+1}+[\mc\otimes_\hc
I^n,\,R]\to\Bbbk$, vanishing on the image of $b:\Omega^1_\hc(R,\,\mc)_\natural/\natural(\mc\otimes_\hc(I^ndR))\to\mc\otimes_\hc R/\mc\otimes_\hc I^{n+1}+[\mc\otimes_\hc
I^n,\,R]$. The image of $b$ coincides with the image of $[,]:\mc\otimes_\hc(R\otimes R)\to\mc\otimes_\hc
R$, projected to $\mc\otimes_\hc R/\mc\otimes_\hc I^{n+1}+[\mc\otimes_\hc
I^n,\,R]$, so one can say that $HC_\hc^{2n}(A,\,\mc)$ is generated by the
linear functionals on $\mc\otimes_\hc R/\mc\otimes_\hc I^{n+1}+[\mc\otimes_\hc
R,\,R]$.

Further observe that two such functionals give the same class in
$HC_\hc^{2n}(A,\,\mc)$, iff their difference is equal to the
composition of $\natural d:\mc\otimes_\hc R/\mc\otimes_\hc I^{n+1}+[\mc\otimes_\hc
R,\,R]\to\Omega^1_\hc(R,\,\mc)_\natural/\natural(\mc\otimes_\hc(I^ndR))$ and
a linear functional $T$ on the space of $1$-forms. Let $\tau_0$ and $\tau_1$ be
two functionals on $\mc\otimes_\hc R/\mc\otimes_\hc I^{n+1}+[\mc\otimes_\hc
R,\,R]$ defining the same cohomology class and $\tilde\tau_0,\ \tilde\tau_1$ be
the functionals on $\mc\otimes_\hc(R/I^{n+1})=\mc\otimes_\hc R/\mc\otimes_\hc I^{n+1}$,
equal to the composition of $\tau_0,\ \tau_1$ with the natural projection $\mc\otimes_\hc R/\mc\otimes_\hc I^{n+1}\to
\mc\otimes_\hc R/\mc\otimes_\hc I^{n+1}+[\mc\otimes_\hc R,\,R]$. Similarly
let $\tilde T$ be the functional on $\Omega^1_\hc(R,\,\mc)/\natural(\mc\otimes_\hc I^ndR)^{-1}$
(here $\natural(\mc\otimes_\hc I^ndR)^{-1}$ denotes the preimage of
$\natural(\mc\otimes_\hc I^ndR)$ in $\Omega^1_\hc(R,\,\mc)$) equal to the
composition of $T$ with the natural projection
$\Omega^1_\hc(R,\,\mc)/(\mc\otimes_\hc
I^ndR)^{-1}\to\Omega^1_\hc(R,\,\mc)_\natural/\natural(\mc\otimes_\hc(I^ndR))$.
Then the condition that $T$ makes $\tau_0$ and $\tau_1$ cohomologous is equivalent
to the equation $\tilde\tau_1-\tilde\tau_0=\tilde Td$.

There is a nice way to interprete this condition: consider the semi-direct
product $L=R\oplus\Omega^1(R)$ of $R$ and $\Omega^1(R)$ (the multiplication is given by the formula
$$ (a,\omega)\cdot(b,\omega') = (ab,a\omega'+\omega b)).$$
Then $L$ is a \hc\/-module
algebra, $\mc\otimes L$ is a
$L$\/-bimodule and (left) \hc\/-module and $I^{n+1}\oplus(I^ndR+I^{n-1}dRI+\dots+dRI^n)$
is an ideal
in $L$ (it is the $n+1$ power of $J=I\oplus\Omega^1(R)$). Consider the 1-parameter family of
homomorphisms of vector spaces $1\otimes_\hc u_t:\mc\otimes_\hc R\to\mc\otimes_\hc L$, extending
the 1-parameter family of homomorphisms of algebras in $\ac^\hc$ (i.e. commuting with the action of \hc) $u_t:R\to L,\ u_t(r)=(r\oplus tdr),\
t\in[0,\,1]$. One can easily check that $u_t(I^{n+1})\subseteq J^{n+1}$.
Define linear functional $\tilde{\tilde T}:\mc\otimes_\hc L/J^{n+1}$ by the formula $\tilde{\tilde T}(m\otimes_\hc(r\oplus xdy))=\tilde\tau_0(m\otimes_\hc r)+\tilde T(m\otimes_\hc xdy),
m\in\mc,\ r\in R,\ xdy\in\Omega^1(R)$. Then clearly $\tilde{\tilde T}$ vanishes
on the space of commutators $[\mc\otimes_\hc L,\,L]$ and on $J^{n+1}$ and verifies the
formulas
$\tilde{\tilde T}u_0=\tilde\tau_0$ and
$\tilde{\tilde T}u_1=\tilde\tau_1$.

Contrarywise, if there exists a 1-parameter family of
morphisms $v_t:(R,\,I)\to(R',\,I')$ in $\ac^\hc$ for a \hc\/-module algebra $R'$
and an \hc\/-module
ideal $I'$ in $R'$, and if there is a linear functional $T'$, vanishing
on $(I')^{n+1}$ and $[\mc\otimes_\hc R',\,R']$, such that $T'v_0=\tilde\tau_0,\
T'v_1=\tilde\tau_1$ (here we write $v_t$ instead of $1\otimes_\hc v_t:\mc\otimes_\hc R\to\mc\otimes_\hc
R'$), then the formula
$$
T(m\otimes_\hc xdy)=\int_0^1T'(m\otimes_\hc(v_t(x)\dot{v}_t(y)))
$$
defines a linear functional on $\Omega^1_\hc(R,\,\mc)$, which vanishes on
the preimage of $\natural(\mc\otimes_\hc I^ndR)$ and on the commutators
$[\mc\otimes_\hc\Omega^1(R),\,R]$ and
such that $\tilde\tau_1-\tilde\tau_0= Td$.

Summarizing these observations, we define {\em \mc\/-twisted \hc\/-equivariant higher even traces of order $n$ on
$A$\/} (or {\em even \mc\/-traces\/} for short) as
\mc\/-twisted traces on $R$-bimodule $\mc\otimes_\hc(R/I^{n+1})$ (defined for arbitrary extension \eqref{eqseq} of
$A$). Here we say that a linear functional on \hc\/-bimodule $M$ over an \hc\/-algebra $R$
is an \mc\/-twisted trace, if it
vanishes on the space of commutators $[\mc\otimes_\hc M,\,R]$ (see formula
\eqref{commute}). We say that two such traces $\tau$ and $\tau'$ defined for
different extensions
$(R,\,I)$ and $(R',\,I')$ are {\em equivalent\/}, if
there exist a third extension $(R'',\,I'')$ and maps $f:R''\to R,\ f':R''\to
R'$ of extensions, such that $\tau f=\tau' f'$.
Since $(RA,IA)$ is a universal extension, we can assume that (mod the equivalence relation)
all traces are functionals on the same linear space.
We say that two traces (on the same extension $(R,\,I)$) are
 {\em homotopic\/}, if
there exists another extension $(R',\,I')$ of $A$ (in $\ac^\hc$), a degree
$n$ even \mc\/-trace $T'$ on
it and a 1-parameter family of homomorphisms $v_t:(R,\,I)\to(R',\,I')$, such
that $T'v_0=\tilde\tau_0,\ T'v_1=\tilde\tau_1$. Then
\begin{prop}
There is a 1-1 correspondence between the cohomology classes in
$HC^{2n}_\hc(A,\,\mc)$ and the homotopy classes of even \mc\/-traces of
degree $n$ on $A$.
\end{prop}

Similarly the degree $2n+1$ cohomology of $A$ can be described in terms of
the odd part of the cohomology of $\xc_\hc^{2n+1}(R,\,I;\,\mc)$ for a
quasi-free extension of $A$. This description would involve the equivalence relation cyclic
1-cocycles on \hc\/-module algebra $R/I^{n+1}$ with coefficients in \mc.
However it is more convenient to use the following ideas similar to the
discussion following the Proposition 9.5 of \cite{CQ2}.

Let $(R,\,I)$ be an extension of $A$ in $\ac^\hc$.
One calls {\em \mc\/-twisted \hc\/-equivariant odd higher trace of degree $n$ on
$A$\/} an $I$\/-adic \mc\/-twisted trace on $\mc\otimes_\hc I^{n+1}$, i.e. a linear
functional on $\mc\otimes_\hc I^{n+1}$ vanishing on the space of
commutators $[\mc\otimes_\hc I^{n+1},\,I]$ (we shall denote the factor-space $\mc\otimes_\hc I^{n+1}/[\mc\otimes_\hc I^{n+1},\,I]$ by $(\mc\otimes_\hc I^{n+1})_\natural$). Two odd traces $\tau$, on $(R,\,I)$ and $\tau'$ on $(R',\,I')$ are equivalent,
if they coincide on an extension $(R'',\,I'')$, which maps to $(R,\,I)$ and
$(R',\,I')$. Two odd traces on the same extension $(R,\,I)$ are {\em homotopic\/}, if
their difference is equal to the restriction of an ordinary trace on $R$.
Then
\begin{prop}
There is a 1-1 correspondence between the the cohomology classes in $H_\hc^{2n+1}(A,\,\mc)$
and homotopy classes of equivalence classes of odd \mc\/-traces of degree $n$ on $A$.
\end{prop}

\subsection{Crainic-type pairing}
\label{sct4}
Now we can define the pairing. Let \hc\/-module algebra $A$ be at the same time
a $C$\/-module algebra in a manner described in section 1, equations \eqref{eqch1}, \eqref{eqch2}.
Let $(R,\,I)$ be an extension of $A$, in
which $R$ is an \hc\/-module and $C$\/-module algebra and $I$ --- ideal in $R$,
stable under the action of $C$ and \hc. We assume that the action of \hc\
and $C$ on $R$ verifies the same conditions \eqref{eqch1} and \eqref{eqch2} and all the
maps in the sequence \eqref{eqseq} are morphisms of \hc\/- and $C$\/-modules.
We also suppose, that the exact sequence \eqref{eqseq} is splittable not
only as the sequence of \hc\/-modules, but also as the sequence of
$C$\/-modules. That is there exists a section $\rho$ of the epimorphism $R\to
A$, such that $\rho(h(a))=h\rho(a)$ and $\rho(c(a))=c\rho(a)$ for all $h\in\hc, c\in C$. One should
think of the universal extension $R=RA,\ I=IA$ of $A$ as the model example
of such extensions.

First of all we interprete the map $\rho$ as a map from $C$ to ${\rm
Hom}\,(B(A),\,R)$. Here we let an element $c\in C$
go to a map $\rho(c):A\to R,\ \rho(c)(a)=\rho(c(a))$ and $B(A)$
stands for the bar-resolution of the algebra $A$ (recall
that as a linear space $B(A)=\displaystyle{\bigoplus_{n\ge1}^\infty}{A^{\otimes n}}$).
The space $\Hom(B(A),\,R)$
is a DG algebra w.r.t. the differential induced from the standard codifferential on
$B(A)$ and the cup-product of maps from coalgebra
$B(A)$ to the algebra $R$: for $\varphi\in\Hom(A^{\otimes p},\,R),\
\psi\in\Hom(A^{\otimes q},\,R)$
we put
\begin{equation}
\label{cupprod}
\varphi\cup\psi:A^{\otimes p+q}\to R,\ \varphi\cup\psi(a_1,\dots, a_{p+q})=\varphi(a_1,\dots, a_p)\psi(a_{p+1},\dots,
a_{p+q}).
\end{equation}
So we can in a unique way extend the map $\rho$ to the map of DG algebras $\rho^\sharp:W(C)\to{\rm
Hom}\,(B(A),\,R)$ (recall that this is the universal property of the Weil algebra, see
remark \ref{rem11}). We begin the description of our pairing construction
with the case of even cohomology classes in $H_\hc^*(A,\,\mc)$.

It is easy to see that the map $\rho^\sharp$ sends elements of the form $\omega_c\in
W(C)$ to $I$ and intertwines the action of \hc\ on $W(C)$ and $\Hom(B(A),\,R)$, if we
define the latter via the action of \hc\/ on $R$. So it gives the maps $\mc\otimes_\hc\rho^\sharp_n:\mc\otimes_\hc W_n(C)\to\mc\otimes_\hc\Hom(B(A),\,R/I^{n+1})\subseteq\Hom(B(A),\,\mc\otimes_\hc
R/I^{n+1})$ and $\mc\otimes_\hc\rho^\sharp_I:I_{n+1}^\hc(C,\,\mc)\to\Hom(B(A),\,\mc\otimes_\hc
I^{n+1})$. Finally, observe that $\mc\otimes_\hc\rho^\sharp$ and $\mc\otimes_\hc\rho^\sharp_I$
send the
commutators' subspace $[\mc\otimes_\hc W(C),\,W(C)]$ and $[\mc\otimes_\hc I(C)^{n},\,I(C)]$
to the corresponding commutators subspace of
morphisms in $\mc\otimes_\hc\Hom(B(A),\,R)$ and $\Hom(B(A),\,R)$, $[\mc\otimes_\hc\Hom(B(A),\,R),\,\Hom(B(A),\,R)]\subseteq\Hom(B(A)^\natural,\,[\mc\otimes_\hc
R,\,R])$ and similarly with the ideal $I$. Here $B(A)^\natural$ denotes the cocommutator subspace
($\Leftrightarrow$ subspace of cyclically invariant tensors) in $B(A)$. Thus we obtain maps of (co)chain complexes
\begin{align}
\tilde\rho_n:W_n(C,\,\mc)_\natural&\to\Hom(B(A)^\natural,\,(\mc\otimes_\hc
R/I^{n+1})_\natural)\\
\intertext{and}
\tilde\rho_{n,I}:(I_{n+1}^\hc(C,\,\mc))_\natural&\to\Hom(B(A)^\natural,\,(\mc\otimes_\hc
I^{n+1})_\natural).
\end{align}
Let $\tau$ be an even \mc\/-trace on $A$ of degree $n$, then $[\tau]$ will denote
the corresponding
cohomology class in $H^{2n}_\hc(A,\,\mc)$. We shall define the product of a class
$[\omega]\in H^p(W^\hc_n(C,\,\mc)_\natural)=H^{p-2n-1}_\hc(C,\,\mc)$ represented by
cocycle $\omega\in W^\hc_n(C,\,\mc)^p_\natural$ with $[\tau]$ by the following rule
\begin{equation}
\label{cupeven}
[\omega]\sharp'[\tau]=[\tau\circ\tilde\rho_n(\omega)\circ N].
\end{equation}
Here $N:B(A)\to B(A)^\natural,\ N=1+t+\dots+t^p$ --- the universal cotrace on $A^{\otimes p}
\subset B(A),\ N:A^{\otimes p}\to (A^{\otimes p})^\natural$.
Thus, the map
on the right can be regarded as cyclic cochain of degree $p-1$ on $A$ (the degree 1 shift
is due to the fact, that the grading in cyclic complex is defined so that degree $k$ is given to linear functionals on $A^{\otimes k+1}$). It is closed,
because $\omega$ is (and $\tilde\rho_n$ commutes with differentials).
\begin{prop}
The operation \eqref{cupeven} is well-defined on the level of cohomology, i.e. it does
not depend on the choice of representative $\omega$ of class $[\omega]$ in
$H^p(W^\hc_n(C,\,\mc)_\natural)$, of trace $\tau$ representing an element
in $H^{2n}_\hc(A,\,\mc)$ and of the splitting $\rho$ of extension
\eqref{eqseq}.
\end{prop}
\begin{proof}
As far as the independence of the class of
$[\tau\circ\tilde\rho_n(\omega)\circ\natural]$ on the choice of $\omega$ is
concerned, it follows from the fact that $\tilde\rho_n$ is a map of chain
complexes.

If $\tau_0\sim\tau_1$, then $\tau_0$ and $\tau_1$ are homotopic; choose a
polynomial 1-parameter family of homomorphisms $u_t:R\to L$ and a degree $n$ \mc\/-trace $T$ on
$L$, connecting $\tau_1$ and $\tau_1$. Then one can consider instead of $u_t$ a
map $U:R\to L[t]$. In this way we obtain from the splitting $\rho$ a map $W(C)\to\Hom(B(A),\,\Omega(\mathbb R^1)\otimes
L)$. Extending this map in a manner similar to what we had above, and using
the trace $\int_0^1\otimes T$ on $\Omega^1(\mathbb R^1)\otimes\mc\otimes_\hc
L/J^{n+1}$ we obtain a homotopy between
$[\tau_0\circ\tilde\rho_n(\omega)\circ N]$ and
$[\tau_1\circ\tilde\rho_n(\omega)\circ N]$.

Similarly, if $\rho^0$ and $\rho^1$ are two splittings of the
sequence \eqref{eqseq}, we combine them in the following way
$\rho=t\rho^0+(1-t)\rho^1,\ t\in[0;\,1]$. As before we use this map
to obtain a homomorphism $W(C)\to\Hom(B(A),\,\Omega(\mathbb
R^1)\otimes R)$, and use the trace $\int_0^1\otimes\tau$ to define
an element connecting $[\tau\circ\tilde\rho^0_n(\omega)\circ N]$ and
$[\tau\circ\tilde\rho^1_n(\omega)\circ N]$.
\end{proof}

Further let $\tau'$ be an odd \mc\/-trace of degree $n$, and
$[\tau']\in H^{2n+1}_\hc(A,\,\mc)$ the cohomology class. We define
its product with the class $[\omega]\in
H^p(I_{n+1}^\hc(C,\,\mc)_\natural)=H^{p-2n-2}_\hc(C,\,\mc)$ (see
remark, following the theorem \ref{theoa}) as the result of the
composition, similar to \eqref{cupeven}:
\begin{equation}
\label{cupodd}
[\omega]\sharp'[\tau']=[\tau'\circ\tilde\rho_{n,I}(\omega)\circ N].
\end{equation}
\begin{prop}
The operation \eqref{cupodd} is well-defined on the level of cohomology.
\end{prop}
\begin{proof}
The only thing that needs proof is the independence of
\eqref{cupodd} of the choice of representative $\tau'$ in the class
$[\tau']$ (independence of the choice of splitting $\rho$ can be
proven similarly to the even case). But $\tau'$ and $\tau''$ are
cohomologous iff their difference is equal to the restriction of a
\mc\/-trace $T$ on $R$. Using this trace instead of $\tau'$ and
$\tau''$, we can define linear functional
$T\circ\tilde\rho_{n,I}(\omega)\circ N$ on $B_p(A)^\natural$, equal
to the difference $\tau'\circ\tilde\rho_{n,I}(\omega)\circ
N-\tau''\circ\tilde\rho_{n,I}(\omega)\circ N$. But $T$ being a trace
on the greater algebra $R$, we can define
$T\circ\tilde\rho_{n,I}(\omega)\circ N$ for arbitrary $\omega'\in
W^\hc(C,\,\mc)_\natural$, in particular we can replace
$I^\hc_{n+1}(C,\,\mc)_\natural$ with $\tilde
I^\hc_{n+1}(C,\,\mc)_\natural$ (see remark following the proof of
theorem \ref{theoa}) and think of the original element $\omega$ as
of a cocycle in $\tilde I^\hc_{n+1}(C,\,\mc)_\natural\subseteq
W^\hc(C,\,\mc)_\natural$. The complex $W^\hc(C,\,\mc)_\natural$
being acyclic, there is an element $\alpha\in
W^\hc(C,\,\mc)_\natural$ for which $\partial\alpha=\omega$. Then the
functional $T\circ\tilde\rho_{n,I}(\alpha)\circ N$ satisfies the
equation $(T\circ\tilde\rho_{n,I}(\alpha)\circ
N)d=\tau'\circ\tilde\rho_{n,I}(\omega)\circ N-
\tau''\circ\tilde\rho_{n,I}(\omega)\circ N$.
\end{proof}

Combining this construction and the result of the previous section we obtain
the desired map:
\begin{align}
\notag
HC^p_\hc(C,\,\mc)\otimes
HC^q_\hc&(A,\,\mc)\stackrel{\sharp'}{\lar}HC^{p+q}(A),\\
\intertext{defined for even $q=2n$ by the formula}
\label{prodeven}
[x]\otimes[y]&\mapsto\tau(\rho_n(\alpha_n^{-1}(x))(Nz))\\
\intertext{and for odd $q=2n+1$}
\label{prododd}
[x]\otimes[y]&\mapsto\tau'(\rho_{n,I}(\beta_n^{-1}(x))(Nz)),
\end{align}
where $\tau$, (resp. $\tau'$) are the even (resp. odd) \mc\/-traces, representing
the class $[y]\in HC^q_\hc(A,\,\mc)$ and $\alpha_n,\ \beta_n$ are the isomorphisms
$H^*(W^\hc_{n+1}(C,\,\mc)_\natural)\cong
HC^{*-1-2n}_\hc(C,\,\mc)$ and
$H^*(I^\hc_{n+1}(C,\,\mc)_\natural)\cong
HC^{*-2-2n}_\hc(C,\,\mc)$ from theorem \ref{theoa} and remark \ref{rem1}, and $z$
is the argument --- element from $B(A)$.

 Our next purpose is to determine the relation of the {\em cup-product\/}
$\sharp'$ of equations \eqref{prodeven}-\eqref{prododd} with
$S$\/-operation on cyclic cohomologies involved. To this end we
shall again generalize the methods, used in \cite{Crain}.

Let once again $\rho$ be a $C$\/- and \hc\/-linear splitting of the
extension \eqref{eqseq}. We consider it as a map from $C$ to
$\Hom(B(A),\,R)$ and extend to a homomorphism of DG algebras
$\rho^\sharp:W(C)\to\Hom(B(A),\,R)$. Consider the universal modules
of 1-forms $\Omega^1(W(C)),\ \Omega^1(R)$ and
$\Omega^1(\Hom(B(A),\,R))$ of algebras $W(C),\ R$ and
$\Hom(B(A),\,R)$ and the universal comodule of 1-forms
$\Omega_1(B(A))$ on coalgebra $B(A)$ (see \cite{Quil2}). First of
all, we extend $\rho^\sharp$ to the map
$\rho_1^\sharp:\Omega^1(W(C))\to\Omega^1(\Hom(B(A),\,R))$ commuting
with differentials and such that the following diagramm commutes
\begin{equation*}
\begin{CD}
{W(C)}           @>{\rho^\sharp}>>     {\Hom(B(A),\,R)}\\
@V{d_u}VV                                @V{d_u}VV\\
{\Omega^1(W(C))}@>{\rho_1^\sharp}>>{\Omega^1(\Hom(B(A),\,R))}.
\end{CD}
\end{equation*}
Here $d_u$ denotes the universal differential for algebras $W(C)$
and $\Hom(B(A),\,R)$. On the other hand, composition of a map
$\phi\in\Hom(B(A),\,R)$ with the universal differential
$d_u:R\to\Omega^1(R)$ and with universal codifferential
$d^u:\Omega_1(B(A))\to B(A)$ are derivatives on $\Hom(B(A),\,R)$
with values in $\Hom(B(A),\,\Omega^1(R))$ and
$\Hom(\Omega_1(B(A)),\,R)$ (observe that both these spaces are
modules over $\Hom(B(A),\,R)$). By the universal property of
$\Omega^1$ there's a map of $\Hom(B(A),\,R)$\/-modules
$\xi:\Omega^1(\Hom(B(A),\,R))\to\Hom(B(A),\,\Omega^1(R))\oplus\Hom(\Omega_1(B(A)),\,R)$
such that $\xi\circ d_u$ is a derivative on $\Hom(B(A),\,R)$.

Let us restrict the image of the map $\xi$ to the first summand and
consider its composition with $\rho_1^\sharp$. It is easy to check
that the resulting map descends to $\rho^\sharp_-
=(\xi\circ\rho_1^\sharp)_\natural:\Omega^1(W(C))_\natural\to\Hom(\Omega_1(B(A))^\natural,\,R_\natural)$
(in other words that commutators in $\Omega^1(W(C))$ send
cocommutators in $\Omega_1(B(A))$ to commutators in $R$). Moreover
if we put
$\rho^\sharp_+=\natural\circ(\rho^\sharp):W(C)\to\Hom(B(A),\,R_\natural)$
then one can check by a straightforward computation that
$X(\rho^\sharp)=(\rho^\sharp_+,\,\rho^\sharp_-)$ is a map of
super-complexes of chain complexes
$X(\rho^\sharp):X(W(C))\to\Hom(X(B(A)),\,R_\natural)$.

Further consider the odd part of the filtration in $X(W(C))$ induced
ideal $I(C)$. The observation that elements of the type $w_c$ in
$W(C)$ are sent by $\rho^\sharp$ to $\Hom(B(A),\,I)$ shows that
$X(\rho^\sharp)\bigl(F^{2n+1}_{I(C)}X(W(C))\bigr)\subseteq\Hom(X(B(A)),\,I^{n+1})$.
Thus we obtain a collection of maps
$X^{2n+1}(\rho^\sharp):X^{2n+1}(W(C),\,I(C))\to\Hom(X(B(A)),\,(R/I^{n+1})_\natural)$.
Finally one can introduce into this construction the coefficients
module \mc, so that all the commutators become \mc\/-twisted
commutators. We obtain a map
\begin{equation}
\label{xrho1}
\tilde X^{2n+1}_\hc(\rho^\sharp):\xc_\hc^{2n+1}(W(C),\,I(C);\,\mc)\to\Hom(X(B(A)),\,(\mc\otimes_\hc
R/I^{n+1})_\natural).
\end{equation}
On the other hand one has the following
\begin{prop}The total (bigraded) complex of $X^{2n+1}_\hc(W(C),\,I(C);\,\mc)$ is
isomorphic to that of the sequence \eqref{comw}.
\end{prop}
The proof of this statement is a word-to-word repetition of Theorem 7.9 of \cite{Crain}), so we omit it.

But the $S$\/-operations in $HC^*_\hc(C,\,\mc)$ are given by diagram
chasing in the sequence \eqref{comw} and similarly the
$S$\/-operations in $HC^*_\hc(A,\,\mc)$ are related to diagram
chasing in the total complex of $X(B(A))$. Thus we see that, if
$q=\mathrm{deg}\,[y]$ is even, one has the equation
$S([x]\sharp'[y])=S[x]\sharp'[y]$. The same is true for odd $q$. To
see this, just observe that instead of dividing out by an ideal we
can restrict to it thus obtaining a map from the exact sequence
\eqref{comi}.

In order to find the relation with $S$\/-operations in
$HC^*_\hc(A,\,\mc)$ we recall that, if $\tau:(\mc\otimes_\hc
R/I^{n+1})_\sharp\to k$ represents $[y]\in HC^{2n}_\hc(A,\,\,\mc)$,
then the class $S[y]$ is represented by the composition
$\tau^p=\tau\circ p_{n+1}$ of $\tau$ with projection
$p_{n+1}:(\mc\otimes_\hc R/I^{n+2})_\sharp\to(\mc\otimes_\hc
R/I^{n+1})_\sharp$. Thus we have $([x]\sharp'
S[y])(b)=\tau^p(\rho_{n+1}(\alpha_{n+1}^{-1}(x))(b))$ (here $b\in
B(A)^\natural$ is the argument of the cyclic cochain). But obviously
$p_{n+1}(\rho_{n+1}(\omega))=\rho_n(\pi_{n+1}(\omega))$ for all
$\omega\in W_{n+1}^\hc(C,\,\mc)_\natural$
($\pi_{n+1}:W_{n+1}^\hc(C,\,\mc)_\natural\to
W_{n}^\hc(C,\,\mc)_\natural$ is the natural projection). So from the
second statement of the theorem \ref{theoa} we obtain
$\tau^p(\rho_{n+1}(\alpha_{n+1}^{-1}(x))(b))=\tau(\rho_n(\pi_{n+1}\circ\alpha_{n+1}^{-1}(x))(b))=
\tau(\rho_n(\alpha_n^{-1}(Sx))(b))$, and hence $[x]\sharp'
S[y]=S[x]\sharp'[y]=S([x]\sharp'[y])$ (the last equality was proven
earlier). Similarly one obtains the same equalities for the odd
case.

We sum up the results of this section in the following theorem
\begin{theorem}
There is a pairing
$$
HC^p_\hc(C,\,\mc)\otimes
HC^q_\hc(A,\,\mc)\stackrel{\sharp'}{\lar}HC^{p+q}(A)
$$
of the Hopf-type cyclic cohomology of \hc\/-module coalgebra $C$ and $C$\/-
and \hc\/-module algebra $A$ with coefficients in a SAYD-module \mc. This
pairing satisfies the following relation with $S$\/-operation
$$
S([x]\sharp'[y])=S[x]\sharp'[y]=[x]\sharp'S[y].
$$
\end{theorem}

\section{Relation with the bivariant theory}
In this section we give a nice interpretation of the construction
described above in terms of a cohomology theory, closely related
with the bivariant cyclic cohomology (see \cite{JK}, \cite{CQ2} and
\cite{CQ3}). In effect, we manage to construct a map from the
cohomology $HC_\hc^*(C,\,\mc)$ to the bivariant cohomology
$HC^*(\xc_A,\,\xc_\hc(R,\,I;\,\mc))$ of two special towers (see \S2
of \cite{CQ2} for the definitions and relation of this cohomology
with that of \cite{JK}). Then we show that the cup-product of
section 3 is equal to the composition of this map with the usual
composition product in bivariant theory. We also use this result to
construct a dual cup-product of cyclic theories.

\subsection{The mapping complex}
For the most part of this section we assume that the algebra $A$ is
finite-dimensional (as $\Bbbk$\/-vector space). First of all observe
that under this condition $\Hom(A,\,V)\cong V\otimes\Hom(A,\,\Bbbk)$
for any vector space $V$. The isomorphism is given by the map
\begin{equation}
\label{isoo} F:V\otimes\Hom(A,\,\Bbbk)\to\Hom(A,\,V),\ F(v\otimes
\varphi)(a)=\varphi(a)v,\ v\in V,\ \varphi\in\Hom(A,\,\Bbbk),\ a\in
A.
\end{equation}
The same is true if we replace $A$ with its bar-resolution,
i.e.
\begin{equation}
\label{isohom} \Hom(B(A),\,V)\cong V\otimes\Hom(B(A),\,\Bbbk)
\end{equation}
as {\em differential graded\/} vector spaces. The differential on
both sides is induced from the differential in the bar-resolution.
To see this recall that
$B(A)=\displaystyle{\bigoplus_{n\ge1}}A^{\otimes n}$ and every
tensor product in this sum is a finite-dimensional space. Let us
denote the chain complex $\Hom(B(A),\,\Bbbk)$ by $B$. Then $B$ is in
fact a DG {\em algebra\/}, multiplication being defined as the
cup-product of cochains for $\varphi\in\Hom(A^{\otimes p},\,\Bbbk),\
\psi\in\Hom(A^{\otimes q},\,\Bbbk)$ we put
$$
\varphi\cup\psi:A^{\otimes p+q}\to \Bbbk,\
\varphi\cup\psi(a_1,\dots, a_{p+q})=\varphi(a_1,\dots,
a_p)\psi(a_{p+1},\dots, a_{p+q})
$$
(compare with \eqref{cupprod}). If $V=R$ is an algebra then one can
check that the isomorphism \eqref{isohom} is in fact an isomorphism
of DG algebras, it is only necessary to show that the extension of
map $F$ in \eqref{isoo} commutes with multiplication, which is
clear.

Now consider the map $\rho^\sharp$ from the previous section (recall
that $\rho^\sharp$ is induced by a splitting $\rho$ of the exact
sequence \eqref{eqseq}). In the view of the previous discussion it
can be interpreted as a homomorphism of the free DG algebra $W(C)$
to $B\otimes R$. In \S14 of \cite{CQ2} it is shown that in case of
usual (non-differential) algebras this homomorphism can be extended
to a unique (up to a homotopy) homomorphism of $X$\/-complexes:
$X(\rho^\sharp):X(W(C))\to X(R)\otimes X(B)$. The following result
is a direct generalization of the corresponding theorems in the
cited paper.
\begin{prop}
\label{prneww} The map $X(\rho^\sharp)$ defined by the methods of
\S14 \cite{CQ2} is a map of super-complexes of cochain complexes.
This map sends the $p$\/-th term of filtration \eqref{filtr} on
$X(W(C))$, associated to the ideal $I(C)\subset W(C)$, to the
subcomplex $X(B)\otimes F_I^{p-1}X(R)$. If all the algebras are
\hc\/-module algebras, then the map $X(\rho^\sharp)$ can be raised
to the homomorphism of \hc\/-equivariant $X$\/-complexes with
coefficients: $X_\hc(\rho^\sharp;\,\mc):X_\hc(W(C),\,\mc)\to
X_\hc(R,\,\mc)\otimes X(B)$.
\end{prop}
\begin{proof}
In order to prove this statement we first recall the construction of
$X(\rho^\sharp)$ from \cite{CQ2}.

First of all one considers the free product $F=R\ast B$ of the
algebras $B$ and $R$. Let $J$ be the ideal in $F$ generated by the
commutators $[b,\,r],\ b\in B,\ r\in R$. Define the algebra $R\,\#B$
as the factor of $F$ by the square of the ideal $J$. Then $R\,\#B$
is a square-zero extension of $R\otimes B\cong F/J$ (we shall denote
the image of $J$ in $R\,\#B$ by the same letter $J$). Since $W(C)$
is a quasi-free (in fact, free) algebra, there's a lifting of
$\rho^\sharp$ to a homomorphism $\rho^\#:W(C)\to R\,\#B$ and one can
use this lifting to define a map of super-complexes $X(W(C))\to
X(R\,\#B)$.

On the other hand, there's an isomorphism of super-complexes
$\xc^2(F,\,J)\cong X(R)\otimes X(B)$. Thus we obtain a morphism:
$$
X(W(C))\to X(R\,\#B)=X(F/J^2)\to\xc^2(F,\,J)\cong X(R)\otimes X(B).
$$
The second map in this composition is the natural projections $F/J^2\to
F/J^2+[J,\,F]$ and
$\Omega^1(F/J^2)_\natural\to\Omega^1(F)_\natural/\natural(JdF)$.

Now the first statement of our proposition follows from the simple
observation that all the algebras and modules in this construction
can be made differential graded and all morphisms will commute with
differentials. For instance one induces the grading on $F=R*B$ from
that on $B$ and similarly the differential. Moreover since $W(C)$ is
free as a DG algebra, one can choose the map $W(C)\to R\,\#B$ to be
a homomorphism of differential algebras.

In order to prove the second statement, let us denote in spite of a
slight abuse of notation by the same symbol $I$ the ideals generated
by $I\subseteq R$ in $F$ and $R\,\#B$. Let $\hat\rho^\sharp:W(C)\to
R\ast B$ denote the homomorphism, covering $\rho^\sharp:W(C)\to
R\otimes B=F/J$ ($\hat\rho^\sharp$ exists, because $W(C)$ is {\em
free\/} DG algebra). Then $\hat\rho^\sharp(I(C))\subseteq I+J$ and
hence $\rho^\#(I^{n+1}(C))\subseteq
I^nJ+I^{n-1}JI+\dots+IJI^{n-1}+JI^n+I^{n+1}$ since $J^2=0$ in
$R\,\#B$. Similarly $\rho^\#([I^n(C),\,W(C)])\subseteq
[I^n+I^{n-1}J+I^{n-2}JI+\dots+JI^{n-1},\,F]$. In the even degree the
isomorphism
\begin{equation}
\label{isoxc}
\phi_\natural:\xc^2(F,\,J)\cong X(R)\otimes X(B)
\end{equation}
is induced by the following map (see \cite{CQ1}, Prop. 1.4 and
\cite{CQ2}, Prop. 14.1):
\begin{align*}
\phi_0:R\ast B&\to R\otimes B\oplus\Omega^1(R)\otimes\Omega^1(B).\\
\intertext{It is given on the generators $r\in R,\ b\in B$ by the
formulas}
\phi_0(r)&=r\otimes1,\ \phi_0(b)=1\otimes b\\
\intertext{and extended to the whole algebra $F=R\ast B$ in a way to obtain a homomorphism of algebras with respect to the following multiplication on the right}
(\xi_0\otimes\eta_0)\circ(\xi_1\otimes\eta_1)\equiv\xi_0\xi_1&\otimes\eta_0\eta_1+(-1)^{|\xi_1|}\xi_0d\xi_1\otimes
d\eta_0\eta_1(\mathrm{mod}\,\Omega^2(R)\otimes\Omega^2(B)).
\end{align*}
Composing this map with the natural projection
$\Omega^1(R)\otimes\Omega^1(B)\to\Omega^1(R)_\natural\otimes\Omega^1(B)_\natural$
(we denote this composition by $\phi_{0,\natural}$) we obtain the
following inclusions
\begin{align*}
\begin{split}
\phi_{0,\natural}\rho^\#\Bigl((F_{I(C)}^{2n+1}X(W(C)))_+\Bigr)
                           &=\phi_{0,\natural}\bigl(\rho^\#(I^{n+1}(C))\bigr)\\
                           &\subseteq\phi_{0,\natural}(I^nJ+I^{n-1}JI+\dots+IJI^{n-1}+JI^n+I^{n+1})\\
                           &\subseteq I^{n+1}\otimes B\oplus\natural(I^ndR)\otimes\Omega^1(B)_\natural\\
                           &\subseteq (I^{n+1}+[I^n,\,R])\otimes B\oplus\natural(I^ndR)\otimes\Omega^1(B)_\natural\\
                           &=(F^{2n}_IX(R)\otimes X(B))_+
\end{split}
\intertext{since the terms containing $J$ give no input to the first
summand and $\phi_0([r,\,b])=dr\otimes db$. Similarly on the even
terms of filtration we obtain}
\begin{split}
\phi_{0,\natural}\rho^\#\Bigl((F_{I(C)}^{2n}X(W(C)))_+\Bigr)
                           &=\phi_{0,\natural}\bigl(\rho^\#(I^{n+1}(C)+[I^{n}(C),\,W(C)])\bigr)\\
                           &\subseteq\phi_{0,\natural}(I^nJ+I^{n-1}JI+\dots+IJI^{n-1}+JI^n+I^{n+1})\\
                           &\quad+\phi_{0,\natural}([I^{n-1}J+\dots+JI^{n-1}+I^n,\,F])\\
                           &\subseteq (I^{n+1}\otimes B\oplus\natural(I^ndR)\otimes\Omega^1(B)_\natural)\\
                           &\quad+([I^n,\,R]\otimes B\oplus\natural(I^{n}dR)\otimes\Omega^1(B)_\natural),
\end{split}\\
\intertext{since the terms of type $[J,\,F]$ are sent to the commutators subspace $[\Omega^1(R),\,R]\otimes\Omega^1(B)$ by $\phi_0$ (see definition). So}
\begin{split}
\phi_{0,\natural}\rho^\#\Bigl((F_{I(C)}^{2n}X(W(C)))_+\Bigr)
                           &\subseteq I^n\otimes B\oplus\natural(I^{n}dR+I^{n-1}dI)\otimes\Omega^1(B)_\natural\\
                           &=(F^{2n-1}_IX(R)\otimes X(B))_+.
\end{split}
\end{align*}
Similarly the odd part of the isomorphism $\phi_\natural$ (equation
\eqref{isoxc}) is given by the following sequence of isomorphisms
(Prop. 14.1 of \cite{CQ2}): first one has
\begin{align*}
\Omega^1(F)\cong F\otimes_R&\Omega^1(R)\otimes_R F\oplus
F\otimes_B\Omega^1(B)\otimes_B F\\
\intertext{which is induced by the universal properties of $\Omega^1(F)$ and the free product. Further
one tensors this equality on both sides by $F/J$ to obtain}
\Omega^1(F)/(JdF+dFJ)&\cong B\otimes\Omega^1(R)\otimes B\oplus
R\otimes\Omega^1(B)\otimes R\\
\intertext{and finally, we factor-out the commutators obtaining}
(\xc^2(F,\,J))_-=\Omega^1(F)_\natural/\natural(JdF)&\cong\Omega^1(R)_\natural\otimes
B\oplus R\otimes\Omega^1(B)_\natural=(X(R)\otimes X(B))_-.
\end{align*}
We shall denote the resulting homomorphism by $\phi_{1,\natural}$.
Then we have the following chain of inclusions (by abuse of notation we shall denote the
natural extension of $\rho^\#$ to $\Omega^1(W(C))_\natural$ by the same
symbol)
\begin{align*}
\begin{split}
\phi_{1,\natural}\rho^\#\bigl((F^{2n}_{I(C)}X(W(C)))_-\bigr)&=\phi_{1,\natural}\rho^\#(\natural(I^n(C)dW(C)))\\
                                                            &\subseteq\phi_{1,\natural}\bigl(\natural((I^n+I^{n-1}J+\dots+JI^{n-1})dF)\bigr)\\
                                                            &\subseteq\phi_{1,\natural}(\natural(I^ndF))\\
                                                            &\subseteq\natural(I^ndR)\otimes B\oplus I^n\otimes\Omega^1(B)_\natural\\
                                                            &\subseteq\natural(I^{n}dR+I^{n-1}dI)\otimes B\oplus I^{n}\otimes\Omega^1(B)_\natural\\
                                                            &=(F^{2n-1}_IX(R)\otimes X(B))_-.
\end{split}\\
\intertext{Since the odd part of the $2n+1$ term of filtration is
equal to a subspace of all 1-forms, generated by $dI^{n+1}$, we
have}
\begin{split}
\phi_{1,\natural}\rho^\#\bigl((F^{2n+1}_{I(C)}X(W(C)))_-\bigr)&\subseteq\phi_{1,\natural}\bigl(\natural(Fd(I^{n+1}+I^nJ+\dots+JI^n))\bigr)\\
                                                              &\subseteq\phi_{1,\natural}(\natural(I^{n+1}dF+I^ndJ)).
\end{split}\\
\intertext{The first summand gives elements in
$\natural(I^{n+1}dR)\otimes B\oplus
I^{n+1}\otimes\Omega^1(B)_\natural$ and the second one --- in
$\natural(I^ndR)\otimes B\oplus
[I^n,\,R]\otimes\Omega^1(B)_\natural$. The last inclusion follows
from the construction of the map $\phi_{1,\natural}$ and the fact
that $J$ is generated by the commutators $[r,\,b]$. So}
\begin{split}
\phi_{1,\natural}\rho^\#\bigl((F^{2n+1}_{I(C)}X(W(C)))_-\bigr)&\subseteq\phi_{1,\natural}(\natural(I^{n+1}dF+I^ndJ))\\
                                                              &\subseteq(\natural(I^{n+1}dR)\otimes B\oplus I^{n+1}\otimes\Omega^1(B)_\natural)\\
                                                              &\quad+(\natural(I^ndR)\otimes B\oplus[I^n,\,R]\otimes\Omega^1(B)_\natural)\\
                                                              &\subseteq\natural(I^ndR)\otimes B\oplus(I^{n+1}+[I^n,\,R])\otimes\Omega^1(B)_\natural\\
                                                              &=(F^{2n}_IX(R)\otimes X(B))_-.
\end{split}
\end{align*}
Finally in order to prove the last statement, one introduces the
action of \hc\ on $F$ via its action on $R$. Then it is a matter of
direct calculation to show that one can extend all the definitions
to the \hc\/-equivariant complexes with coefficients in \mc. For
instance, one should use the maps $1_\mc\otimes_\hc\rho^\sharp$ and
$1_\mc\otimes_\hc\rho^\#$ instead of $\rho^\sharp$ and $\rho^\#$,
regard the map $\phi_0$ as a map of \hc\/-modules, and modify the
definition of $\phi_{1,\natural}$. Observe that the image of
$F^p_{I(C)}X^\hc(W(C),\,\mc)$ is included in
$F^{p-1}_IX^\hc(R,\,\mc)\otimes X(B)$, just like in the
non-equivariant case.
\end{proof}
\begin{rem}\rm
\label{remmorph} One can show as a further generalization of \S14,
\cite{CQ2} (prop. 14.3) that the map $X(W(C))\to X(R)\otimes X(B)$
defined above is a natural extension of the map $X(W(C))\to
R_\natural\otimes X(B)$, induced by the map $X(W(C))\to X(R\otimes
B)$ and the natural projection $\Omega^1(R\otimes B)\to
R\otimes\Omega^1(B)$. Similarly in the equivariant case our map is
an extension of $X(W(C))\to(\mc\otimes_\hc R)_\natural\otimes X(B)$.
And on the $2n+1$\/-st term of filtration our map is an extension of
$F_{I(C)}^{2n+1}X(W(C))\to I^{n+1}/I^{n+1}\bigcap[R,\,R]\otimes
X(B)$, while on the $2n$\/-th term $F_{I(C)}^{2n}X(W(C))\to
I^{n+1}/[I^n,\,I]\otimes X(B)$, since
$\natural(I^{n+1}dR)\subseteq\natural(I^ndI)$ and hence
$\natural(I^{n+1}dR+I^ndI)=\natural(I^ndI)$.
\end{rem}
From the proposition \ref{prneww} it follows, that one can define
the maps of super-complexes of cochain complexes for all $p\ge1$:
\begin{equation}
\xc^p_\hc\rho_\mc:\xc^p_\hc(W(C),\,I(C);\,\mc)\lar\xc^{p-1}_\hc(R,\,I;\,\mc)\otimes
X(B).
\end{equation}
Since we have assumed that $A$ is finite-dimensional, then similarly
to the discussion in the beginning of this section we have the
isomorphism {\em of towers of complexes\/}
$\xc_\hc(R,\,I;\,\mc)\otimes X(B)$ and
$\Hom(X(B(A)),\,\xc_\hc(R,\,I;\,\mc))$ given by the level-wise
isomorphisms
\begin{equation}
\label{xiso} \xc^{p}_\hc(R,\,I;\,\mc)\otimes
X(B)\cong\Hom(X(B(A)),\,\xc^{p}_\hc(R,\,I;\,\mc)).
\end{equation}
Let us denote the tower of super-complexes $\{\Hom(X(B(A)),\,\xc^{p}_\hc(R,\,I;\,\mc))\}$
by $\xc_\hc(A;(R,\,I;\,\mc))$. We come to the following proposition
\begin{prop}
\label{propsup} There exists a morphism of towers of super-complexes
of cochain complexes
$\widetilde\xc_\hc\rho_\mc:\xc_\hc(W(C),\,I(C);\,\mc)\to\xc_\hc(A;(R,\,I;\,\mc))$
induced by a \hc\/- and $C$\/-linear splitting $\rho$ of the exact
sequence \eqref{eqseq}. By this we mean that there exists a series
of morphisms $\widetilde\xc^p_\hc\rho_\mc,\ p\ge1$,
\begin{equation}
\label{xrho22}
\widetilde\xc^p_\hc\rho_\mc:\xc^p_\hc(W(C),\,I(C);\,\mc)\lar\Hom(X(B(A)),\,\xc^{p-1}_\hc(R,\,I;\,\mc)),\
p\ge1,
\end{equation}
which commute with the natural projections of complexes on both sides:
\begin{equation*}
\begin{CD}
{\xc^p_\hc(W(C),\,I(C);\,\mc)}@>{\widetilde\xc^p_\hc\rho_\mc}>>{\Hom(X(B(A)),\,\xc^{p-1}_\hc(R,\,I;\,\mc))}\\
@VVV                                                                 @VVV\\
{\xc^{p-1}_\hc(W(C),\,I(C);\,\mc)}@>{\widetilde\xc^{p-1}_\hc\rho_\mc}>>{\Hom(X(B(A)),\,\xc^{p-2}_\hc(R,\,I;\,\mc))}.
\end{CD}
\end{equation*}
Moreover the homotopy classes of these morphisms does not depend on
the choice of the splitting $\rho$.
\end{prop}
Before we prove this statement, we need to explain in what way we understand the homotopy between the maps of
super complexes in this setting. If $Y=\Bigl\{(\{Y_0^n\}_{n=0}^\infty,\,d)\raisebox{1pt}{$\begin{smallmatrix}\stackrel{\scriptstyle{b_0}}{\scriptstyle\longrightarrow}\\
\stackrel{\scriptstyle\longleftarrow}{\scriptstyle{
b_1}}\end{smallmatrix}$}(\{Y_1^n\}_{n=0}^\infty,\,d)\Bigr\}$ is a supercomplex
of cochain complexes, then we define its cohomology as the shifted by $-1$ cohomology of the
following periodic double complex
\begin{equation}
\label{supcom}
\begin{CD}
 \dots    @.          \dots   @.         \dots   @.        \dots   @.         \dots   @.    {}   \\
@A{d_0}AA             @A{d_1}AA          @A{d_0}AA         @A{d_1}AA          @A{d_0}AA     @.    \\
{Y_0^2} @>{b_0}>>     {Y_1^2} @>{b_1}>>  {Y_0^2} @>{b_0}>> {Y_1^2} @>{b_1}>>  {Y_0^2} @>>>  \dots\\
@A{d_0}AA             @A{d_1}AA          @A{d_0}AA         @A{d_1}AA          @A{d_0}AA     @.    \\
{Y_0^1} @>{b_0}>>     {Y_1^1} @>{b_1}>>  {Y_0^1} @>{b_0}>> {Y_1^1} @>{b_1}>>  {Y_0^1} @>>>  \dots\\
@A{d_0}AA             @A{d_1}AA          @A{d_0}AA         @A{d_1}AA          @A{d_0}AA     @.    \\
{Y_0^0} @>{b_0}>>     {Y_1^0} @>{b_1}>>  {Y_0^0} @>{b_0}>> {Y_1^0} @>{b_1}>>  {Y_0^0} @>>>  \dots.
\end{CD}
\end{equation}
Clearly any map $X\to Y$ of two super-complexes of cochain complexes
defines a map of bicomplexes \eqref{supcom} and hence of their total
complexes. The proposition \ref{propsup} states that the induced
maps of total complexes are homotopic as maps of cochain complexes.
\begin{proof}
The first and the second statements of this proposition are evident.
(They follow from the inclusions of the terms of filtration, Prop.
\ref{prneww}.) In order to proof that the maps
$\widetilde\xc^p_\hc\rho_\mc$ and $\widetilde\xc^p_\hc\rho'_\mc$
defined for two different splittings $\rho$ and $\rho'$ are
homotopic, we consider as in section \ref{sct4} linear combination
of splittings $t\rho+(1-t)\rho'$ as a map from $C$ to the algebra of
$R$\/-valued polynomials $R[t]$. We extend this map to a
homomorphism $X(W(C))\to\Hom(X(B(A)),\,X(\Omega^*(\mathbb
R^1)\otimes R))$ and similarly in \hc\/-equivariant setting and with
coefficients module \mc. Finally we consider the corresponding maps
of total complexes \eqref{supcom} and integrate the forms on
$\mathbb R^1$ that belong to the image of this map from $0$ to $1$.
This gives the desired homotopy.
\end{proof}

Let us conclude this subsection with a brief discussion of corrections that
should be made in our construction, if the algebra $A$ is not
finite-dimensional.

First of all if $A$ can be represented as an inverse limit of its
finite-dimensional subalgebras $A_\lambda$, then everything is
clear: one can obtain all the maps that we need as direct and
inverse limits of the corresponding maps for the subalgebras
$A_\lambda$. In general $A$ can not be represented in this way.
However its bar-resolution $B(A)$ is a direct limit of its
finite-dimensional differential subcoalgebras $C_\lambda,\
\lambda\in\Lambda$. Moreover one can choose the subcoalgebras
$C_\lambda$ so that $C_\lambda\subseteq C_{\lambda'},\
\lambda\preccurlyeq\lambda'$. Then $B=\Hom(B(A),\,k)$ is an inverse
limit of its finite-dimensional subalgebras
$B_\lambda=\Hom(C_\lambda,\,k)$. For each of these subalgebras
there's a map
$$
\xc_\hc(R,\,\mc)\otimes X(B_\lambda)\to\Hom(X(C_\lambda),\,\xc_\hc(R,\,\mc)).
$$
These maps commute with natural projections $B_\lambda\to B_{\lambda'}$ and
thus define a map
$$
\xc_\hc(R,\,\mc)\otimes\lim_{\lal}X(B_\lambda)\to\lim_{\lal}(\xc_\hc(R,\,\mc)\otimes
X(B_\lambda))\to\lim_{\lal}\Hom(X(C_\lambda),\,\xc_\hc(R,\,\mc)).
$$
Combining the map $W^\hc(C,\,\mc)\to\Hom(B(A),\,R)$ with natural projections $\Hom(B(A),\,R)\to\Hom(C_\lambda,\,R)\cong B_\lambda\otimes
R$ and passing to the \xc\/-complexes we obtain a collection of maps
$$
X_\hc(W(C),\,\mc)\to\xc_\hc(R,\,\mc)\otimes X(B_\lambda),\
\lambda\in\Lambda,
$$
which verify the same properties with respect to the filtration as the map of
Prop. \ref{prneww} and commute with projections in the inverse system. Thus
we obtain a map
$$
X_\hc(W(C),\,\mc)\to\lim_{\lal}\Hom(X(C_\lambda),\,\xc_\hc(R,\,\mc)).
$$
Finally we observe that if $C_\lambda$ are chosen as described above,
there's a map
\begin{equation}
\begin{split}
&\displaystyle{\lim_{\lal}}\Hom(X(C_\lambda),\,\xc_\hc(R,\,\mc))\to\Hom(\displaystyle{\lim_{\lar}}X(C_\lambda),\,\xc_\hc(R,\,\mc))\\
&\to\Hom(X(\displaystyle{\lim_{\lar}}C_\lambda),\,\xc_\hc(R,\,\mc))=\Hom(X(B(A)),\,\xc_\hc(R,\,\mc)).
\end{split}
\end{equation}

\subsection{Towers of supercomplexes and bivariant cohomology}
It is natural to regard the tower $\xc_\hc(A;(R,\,I;\,\mc))$ as complex, determining a
(version of) the bivariant cyclic cohomology theory of algebra $A$ with values in $A$ (with coefficients \mc).
Below we shall develope this analogy and show, how one obtains a map of two
towers of super-complexes from elements of $\xc_\hc(A;(R,\,I;\,\mc))$.
However this construction of bivariant cohomology is somewhat different from those, given in \cite{CQ2} and
\cite{JK}. We are not going to discuss here in details the relation of
these two cohomology theories (that of the definition \ref{defbivth} below and the
one defined in the cited papers).

First observe that $X(B(A))$ is a super-complex of {\em chain\/}
complexes (i.e. $X(B(A))=\bigl\{X_0\leftrightarrows X_1\bigr\}$
where both $X_i,\ i=0,1$ are complexes with differentials of degrees
$-1$). We shall denote the structure maps of this complex by
$\beta:B(A)\to\Omega_1(B(A))^\natural,\
\delta\natural:\Omega_1(B(A))^\natural\to B(A)$ and the symbols
$\partial_0$ and $\partial_1$ will denote the differentials induced
from the standard differential in $B(A)$. One can construct the
$2$\/-periodic bicomplex $biX(B(A))$ in a manner similar to
\eqref{supcom}. Consider the total complex
$Q=\mathrm{Tot}\bigl(biX(B(A))\bigr)$. Now cyclic cohomology
$HC^p(A)$ of $A$ is equal to the $(p-2)$\/-dimensional cohomology of
$Q$ (the first column of $biX(B(A))$, which is just $B(A)$, is
acyclic and the $n$\/-th tensor power of $A$ corresponds to the
$n-1$ degree in cyclic homology).

On the other hand, we can associate to $X(B(A))$ a tower of $2$\/-periodic
bicomplexes $X^p(A),\ p\ge0,\ X^p(A)=\{X^p(A)_{i,j}\},\ i,j\in\mathbb Z,\
j\ge1$, where
\begin{equation}
X^p(A)_{i,j}=\begin{cases}
             0,\ &j>p+1\\
             X_\nu(B(A))_{p+1}/\partial_{\nu}X_\nu(B(A))_{p+2},\ &j=p+1,\ \nu=i+2\mathbb Z\\
             X_\nu(B(A))_j\ &1\le j<p+1,\ \nu=i+2\mathbb Z
            \end{cases}
\end{equation}
(here $\partial_i$ denotes the differential in the $i$\/-th column). For
example $X^1(A)$ looks as follows
$$
\begin{CD}
 {} @.    \dots               @.            \dots                                       @.                     \dots               @.            \dots                                        @.     {} \\
 @.       @VVV                              @VVV                                                               @VVV                              @VVV                                                @. \\
 {} @.  0                   @<<<          0                                           @<<<                   0                   @<<<          0                                            @<<<   {} \\
 @.       @VVV                              @VVV                                                               @VVV                              @VVV                                                @. \\
 {} @.  {B_2/\partial_0B_3} @<{\beta'}<<  {\bar\Omega(B_2)/\partial_1\bar\Omega(B_3)} @<{\delta\natural'}<<  {B_2/\partial_0B_3} @<{\beta'}<<  {\bar\Omega(B_2)/\partial_1\bar\Omega(B_3)}  @<<<   {} \\
 @.       @V{\partial_0'}VV                 @V{\partial_1'}VV                                                  @V{\partial_0'}VV                 @V{\partial_1'}VV                                   @. \\
 {} @.  {B_1}               @<{\beta}<<   {\bar\Omega(B_1)}                           @<{\delta\natural}<<   {B_1}               @<{\beta}<<   {\bar\Omega(B_1)}                            @<<< {.}
\end{CD}
$$
Here we abbreviated $B_p(A)$ to $B_p$ and
$\Omega_1(B_p(A))^\natural$ to $\bar\Omega(B_p)$, $\beta',\
\delta\natural'$ and $\partial_i'$ denote the maps induced by
$\beta',\ \delta\natural'$ and $\partial_i'$ on factor-spaces. The
leftmost column at this diagram corresponds to $i=0$ and the lowest
row is indexed by $j=1$. The maps $X^{p+1}(A)\to X^p(A),\ p\ge0$ are
given by the natural projections of columns.

Observe now that the total complex of $X^p(A)$ is $2$\/-periodic and
we can finally define the tower of super-complexes $\xc(B(A))$ by
the formula
\begin{equation}
\xc^p_{\nu+2\mathbb Z}(B(A))=\mathrm{Tot}_{\nu}(X^p(A)),\ \nu=0,1.
\end{equation}
Here $\mathrm{Tot}_q(X^p(A))=\bigoplus_{i+j=\nu}X^p(A)_{i,j}$. The differentials in $\xc^p(B(A))$ are induced from $\mathrm{Tot}_*(X^p(A))$. Then we have
the following proposition.
\begin{prop}
\label{propxb}
The tower $\xc(B(A))$ is a special tower of super-complexes, whose
cohomology is given by equation
\begin{align*}
H_\nu(\xc^p(B(A)))&=\begin{cases}
                    HC_p(A),\ &\nu=p+2\mathbb Z\\
                    HD_p(A),\ &\nu=p-1+2\mathbb Z
                    \end{cases}\\
\intertext{and the $S$\/-operation in the cyclic cohomology of $A$ is induced by
the natural projections $\xc^{p+2}\to\xc^p$. Further the periodic cohomology
of $A$ is given by the formula}
H_\nu(\widehat\xc(B(A)))&=HP_\nu(A)
\end{align*}
where
$\widehat\xc(B(A))={\displaystyle{\lim_{\longleftarrow}}}\,\xc^p(B(A))$.
In particular the supercomplex $\widehat\xc(B(A)))$ with its natural
filtration is homotopy equivalent to the inverse limit of the tower
$\xc_A$ of Cuntz and Quillen.
\end{prop}
\begin{proof}
In order to prove the first equality we observe that $H_{p+2\mathbb
Z}(\xc^p(B(A)))=H_{p+2}(Q)$, which is equal, as we have remarked, to
$HC_p(A)$. In fact $\xc^p_{p+2\mathbb Z}(B(A))=Q_{p+2}/(B(A)\bigcap
Q_{p+2}+\partial_1\Omega_1(B_{p+2}(A))^\natural)$, and we recall,
that $B(A)$ is acyclic.

In order to prove two other identities we observe that the bicomplex
$X^p(A)$ is the total complex of $X(B(A))/F^p(A)$, where
$F^0(A)\supset F^1(A)\supset\dots\supset F^p(A)\supset\dots\supset
X(B(A))$ is a decreasing filtration,
\begin{equation}
F^p(A)_0=\partial_0B_{p+2}(A)\oplus\bigoplus_{k\ge p+2}B_k(A),\quad
F^p(A)_1=\partial_1\Omega_1(B_{p+2}(A))^\natural\oplus\bigoplus_{k\ge
p+2}\Omega_1(B_k(A))^\natural.
\end{equation}
Now it is evident that the homology of
$\mathrm{gr}^p\xc(B(A))=F^{p-1}(A)/F^p(A)$ is equal to
\begin{equation}
H_{p+2\mathbb Z}(\mathrm{gr}^p\xc(B(A)))=HH_p(A),\qquad H_{p-1+2\mathbb
Z}(\mathrm{gr}^p\xc(B(A)))=0.
\end{equation}
Hence the tower $\xc(B(A))$ is special. This observation also proves the
formulas for the periodic and de Rham cohomology of $\xc(B(A))$.

Finally the statement, concerning the supercomplex
$\widehat\xc(B(A))$ comes from the explicit description of this
complex as the super complex, associated with the bicomplex
$biX(B(A))$. So the complex $\widehat\xc(B(A))$ is equivalent (in
fact, isomorphic) to the standard complex computing the periodic
cyclic homology, which is known to be homotopy equivalent to that of
Cuntz and Quillen.
\end{proof}
\begin{rem}\rm
Explicitly complexes $\xc^p(B(A))$ are given by the following formulas: if $p=2k$
\begin{align}
\begin{split}
\xc^{2k}(B(A))_0&=\Omega_1(B_{2k+1}(A))^\natural/\partial_1\Omega_1(B_{2k+2}(A))^\natural\oplus B_{2k}(A)\\
                &\quad\oplus\Omega_1(B_{2k-1}(A))^\natural\oplus B_{2k-2}(A)\oplus\dots\oplus\Omega_1(B_1(A))^\natural,
\end{split}\\
\begin{split}
\xc^{2k}(B(A))_1&=B_{2k+1}(A)/\partial_0B_{2k+2}(A)\oplus\Omega_1(B_{2k}(A))^\natural\\
                &\quad\oplus B_{2k-1}(A)\oplus\Omega_1(B_{2k-2}(A))^\natural\oplus\dots\oplus B_1(A),
\end{split}\\
\intertext{and if $p=2k+1$}
\begin{split}
\xc^{2k+1}(B(A))_0&=B_{2k+2}(A)/\partial_0B_{2k+3}(A)\oplus\Omega_1(B_{2k+1}(A))^\natural\\
                  &\quad\oplus B_{2k}(A)\oplus\Omega_1(B_{2k-2}(A))^\natural\oplus\dots\oplus\Omega_1(B_1(A))^\natural,
\end{split}\\
\begin{split}
\xc^{2k+1}(B(A))_1&=\Omega_1(B_{2k+2}(A))^\natural/\partial_1\Omega_1(B_{2k+3}(A))^\natural\oplus B_{2k+1}(A)\\
                  &\quad\oplus\Omega_1(B_{2k}(A))^\natural\oplus B_{2k-2}(A)\oplus\dots\oplus B_1(A).
\end{split}
\end{align}
\end{rem}
In a dual way, given a super-complex of {\em cochain complexes\/}
$Y^0\leftrightarrows Y^1$ (i.e. a supercomplex in which both
summands $Y^0$ and $Y^1$ are positively graded complexes
$Y^\nu=\{Y^\nu_p\}_{p\ge1}$ with differentials $\partial^\nu$ of
degree $+1,\ \nu=0,1$, see above) one can associate to it an {\em
inverse tower\/} of supercomplexes:
\begin{align}
\yc_p^\nu&=\mathrm{Tot}^{\nu}\tilde Y_p,\ \nu=0,1\\
\intertext{where $\tilde Y_p,\ p\ge 0$ is a tower of bicomplexes, $\tilde Y_p=\{\tilde Y_p^{i,j}\}$}
\tilde Y_p^{i,j}&=\begin{cases}
             0,\ &j>p\\
             \mathrm{Ker}\partial^\nu Y^\nu_{p},\ &j=p,\ \nu=i+2\mathbb Z\\
             Y^\nu_j,\ &0\le j<p,\ \nu=i+2\mathbb Z.
                  \end{cases}
\end{align}
Here the term {\it inverse tower\/} means that this time the
structure maps act in the positive direction, $\yc^p\to\yc^{p+1}$,
i.e. their direction is opposite to that in \xc\/-tower. Explicitly
the \hbox{$p=2k$\/-th} level of this tower is given by
\begin{align}
\yc_{2k}^0&= \mathrm{Ker}\partial^0(Y^0_{2k})\oplus Y_{2k-1}^1\oplus Y_{2k-2}^0\oplus\dots\oplus Y_0^0,\\
\yc_{2k}^1&= \mathrm{Ker}\partial^1(Y^1_{2k})\oplus Y_{2k-1}^0\oplus Y_{2k-2}^1\oplus\dots\oplus Y_0^1,\\
\intertext{and the $p=(2k+1)$\/-st level by}
\yc_{2k+1}^0&= \mathrm{Ker}\partial^1(Y^1_{2k+1})\oplus Y_{2k}^0\oplus Y_{2k-1}^1\oplus\dots\oplus Y_0^0,\\
\yc_{2k+1}^1&= \mathrm{Ker}\partial^0(Y^0_{2k+1})\oplus Y_{2k}^1\oplus Y_{2k-1}^0\oplus\dots\oplus Y_0^1.\\
\end{align}
In particular one can apply this construction to
$\xc^k_\hc(A;(R,\,I;\,\mc))$. In this case
$Y^0_p=\Hom(B_{p+1}(A),\,\xc^k(R,\,I;\,\mc)_0)\oplus\Hom(\Omega_1(B_{p+1}(A))^\natural,\,\xc^k(R,\,I;\,\mc)_1)$
and
$Y^1_p=\Hom(B_{p+1}(A),\,\xc^k(R,\,I;\,\mc)_1)\oplus\Hom(\Omega_1(B_{p+1}(A))^\natural,\,\xc^k(R,\,I;\,\mc)_0)$.
Since the differentials $\partial^\nu$ in this complex are given by
dualizing the corresponding differentials $\partial_\nu$ in
$X(B(A))$, we conclude that the following statement holds
\begin{prop}
Let $\yc^k_p(A,\,R,\,I;\hc,\,\mc)$ denote the $p$\/-th level of the tower,
associated to $\xc^k_\hc(A;(R,\,I;\,\mc))$. Then there is an isomorphism of
supercomplexes
$$
\yc^k_p(A,\,R,\,I;\hc,\,\mc)\lar\Hom(\xc^p(B(A)),\,\xc^k_\hc(R,\,I;\,\mc)).
$$
Here on the right hand side stands the mapping complex of two
supercomplexes.
\end{prop}
Observe that the collection of supercomplexes
$\yc^k_p(A,\,R,\,I;\hc,\,\mc)$ is naturally a tower in both
directions, i.e. there are maps
$\yc^k_p(A,\,R,\,I;\hc,\,\mc)\to\yc^k_{p+1}(A,\,R,\,I;\hc,\,\mc)$
and
$\yc^k_p(A,\,R,\,I;\hc,\,\mc)\to\yc^{k-1}_p(A,\,R,\,I;\hc,\,\mc)$,
commuting with differentials. Consider the diagram (we shall call it
{\em a bitower, associated to $\yc^k_p(A,\,R,\,I;\hc,\,\mc)$\/})
\begin{equation}
\label{diabiv}
\begin{CD}
{\dots}   @.   {\dots}   @.   {\dots}   @.   {\dots}   @.   {}      \\
@VVV           @VVV           @VVV           @VVV           @.      \\
{\yc^2_0} @>>> {\yc^2_1} @>>> {\yc^2_2} @>>> {\yc^2_3} @>>> {\dots} \\
@VVV           @VVV           @VVV           @VVV           @.      \\
{\yc^1_0} @>>> {\yc^1_1} @>>> {\yc^1_2} @>>> {\yc^1_3} @>>> {\dots} \\
@VVV           @VVV           @VVV           @VVV           @.      \\
{\yc^0_0} @>>> {\yc^0_1} @>>> {\yc^0_2} @>>> {\yc^0_3} @>>> {\dots,}
\end{CD}
\end{equation}
where $\yc_p^k$ stands for $\yc^k_p(A,\,R,\,I;\hc,\,\mc)$. Let $T_n\yc,\ n\in\mathbb Z$
denote the subdiagram of \eqref{diabiv}, such that $\yc_p^k\in T_n\yc$ iff
$p-k\ge n$. Now we can give the following definition
\begin{df}
\label{defbivth}
The $n$\/-dimensional bivariant cyclic cohomology of the tower
$\xc_\hc(A;(R,\,I;\,\mc))$ is
\begin{equation}
\label{eqbivcoh}
HC^n(A,\,R,\,I;\hc,\,\mc)\eqdef H_{n+2\mathbb
Z}(\lim_{\stackrel{\longleftarrow}{k}}\lim_{\stackrel{\lar}{p}}(\yc^k_p\in T_n\yc)).
\end{equation}
\end{df}
This definition is inspired by discussions in \cite{CQ3} and
\cite{CQ2}. Recall (\cite{CQ3}), that for two pro-vector spaces
$V=\{V_i\}$ and $W=\{W_j\}$ one defines the space of homomorphisms
$V\to W$ by the formula
$$
\Hom(V,\,W)=\lim_{\stackrel{\longleftarrow}{j}}\lim_{\stackrel{\lar}{i}}\Hom(V_i,\,W_j).
$$
Further, given two towers of supercomplexes $\xc$ and $\xc'$ (for
instance, one can take the canonical towers associated to two
algebras), one defines their bivariant cohomology (\cite{CQ2}) by
formula
\begin{equation}
HC^n(\xc,\,\xc')\eqdef H_{n+2\mathbb
Z}(\Hom^n(\widehat\xc,\,\widehat\xc')),
\end{equation}
where $\widehat\xc=\displaystyle{\lim_{\longleftarrow}}\xc^p$ and
\begin{equation}
\Hom^n(\widehat\xc,\,\widehat\xc')=\Bigl\{f:\xc\to\xc'|\forall m,\ f(F^{m+n}\widehat\xc)\subset
F^m\widehat\xc'\Bigr\}.
\end{equation}
Here $F^n\widehat\xc$ is the natural filtration on the inverse limit of
towers.

Now the following proposition is an evident consequence of our definitions
and Proposition \ref{propxb} (its last statement):
\begin{prop}
\label{propbiv1}
There's a natural map of complexes
\begin{equation}
\label{eqbiv1}
\lim_{\stackrel{\longleftarrow}{k}}\lim_{\stackrel{\lar}{p}}(\yc^k_p\in
T_n\yc)\to\Hom^n(\widehat\xc(B(A)),\,\widehat\xc_\hc(R,\,I;\,\mc)),
\end{equation}
inducing (since $\widehat\xc(B(A))$ is equivalent to $\widehat\xc_A$) a map of cohomology $HC^*(A,\,R,\,I;\hc,\,\mc)\to
HC^n(\xc_A,\,\xc_\hc(R,\,I;\,\mc))$.
\end{prop}
In fact one can make the following conjecture
\begin{conj}
There's an isomorphism of homology theories
$$
HC^*(A,\,R,\,I;\hc,\,\mc)\cong HC^*(\xc_A,\,\xc_\hc(R,\,I;\,\mc))
$$
induced by the map \eqref{eqbiv1} of Proposition \ref{propbiv1}.
\end{conj}
However we shall not need this statement in full generality.

In fact the construction described above can be applied to any tower
of supercomplexes of cochain complexes. Let $X(k)$ be such a tower,
$k\in\mathbb N$ is the level of the tower and for every $k$,
$X(k)=\{X_\nu^p(k)\},\ p\in\mathbb N,\ \nu=0,1$. To each
supercomplex (of chain complexes) $X(k)$ we associate an inverse
tower of supercomplexes $\yc_p(X(k))=\yc'_p{}^k$ (we write $\yc'$ to
distinguish this case from $\yc^k_p(A,\,R,\,I;\hc,\,\mc)$). Finally
we define the cyclic cohomology of the tower $X(k)$ by the formula
similar to \eqref{eqbivcoh}:
\begin{equation}
\label{cyccohtow}
HC^n(X(k))\eqdef H_{n+2\mathbb
Z}(\lim_{\stackrel{\longleftarrow}{k}}\lim_{\stackrel{\lar}{p}}(\yc^k_p\in T_n\yc)).
\end{equation}

In particular let us apply this construction to the tower
$\xc_\hc(W(C),\,I(C);\,\mc)=\{\xc^k_\hc(W(C),\,I(C);\,\mc)\}$ of
supercomplexes of cochain complexes. Then it is clear that for a
fixed $k$ the cohomology of the \yc\/-tower
$\yc_p(\{\xc^k_\hc(W(C),\,I(C);\,\mc)\})$ associated to
$\xc^k_\hc(W(C),\,I(C);\,\mc)$ coincides with the cohomology of the
bicomplex constructed from $\xc^k_\hc(W(C),\,I(C);\,\mc)$, so that
the injection $\yc_p\to\yc_{p+2}$ gives the $S$\/-operation in this
cohomology. The following proposition can be proven by mere
inspection of definitions (take care of the dimension shift in our
definitions of towers associated to $\xc_\hc(W(C),\,I(C);\,\mc)$ and
$X(B(A))$).
\begin{prop}
\label{propmapcomp1}
The map $\widetilde\xc^p_\hc\rho_\mc$ (formula \eqref{xrho22}) determines a map of towers
\begin{equation}
\label{eqmapcomp1}
\yc_p(\{\xc^k_\hc(W(C),\,I(C);\,\mc)\})\to\yc_{p-1}^{k-1}(A,\,R,\,I;\hc,\,\mc).
\end{equation}
and hence a map of cohomology $H^*(\{\xc^k_\hc(W(C),\,I(C);\,\mc)\})\to H^*(A,\,R,\,I;\hc,\,\mc)$
(on the left here stands the homology of the double tower, associated with $W(C)$,
see the discussion preceeding this proposition).
\end{prop}

Consider now an element $[x]\in HC^m_\hc(C,\,\mc)$. With the help of the isomorphisms $\alpha_n$
(Theorem \ref{theoa}) we can associate to it a series of elements $\alpha_n^{-1}([x])\in H^{m+2n+1}(C,\,\mc;\,n)=
H^{m+2n+1}(\xc_\hc^{2n+1}(W(C),\,I(C);\,\mc))$ and cocyclces $\{\alpha_n^{-1}([x])\}$. From the discussion
of $S$\/-operations in the complexes $W_k^\hc(C,\mc)_\natural$ and supercomplex $\xc^k_\hc(W(C),\,I(C);\,\mc)$
(see propositions \ref{longcom} and \ref{shortcom} and the paragraph above),
we conclude, that the cocycles $\{\alpha_n^{-1}([x])\}$ fit together.
In fact it was shown that on the cochain level the $S$-operator defined via the diagram
(\ref{comw}) (which gives the even horizontal arrows in diagram~(\ref{diabiv})) coincides with
the projections $\pi_n : W_n^\hc(C,M)_\natural\to W_{n-1}^\hc(C,M)_\natural$ (which give
the even vertical arrows in~(\ref{diabiv})). Thus we obtain a
cocycle in the limit complex $\displaystyle{\lim_{\stackrel{\longleftarrow}{k}}}\displaystyle{\lim_{\stackrel{\lar}{p}}}(\yc^k_p\in
T_n\yc)$ (where $\yc$ is the $\yc$\/-bitower associated to
$\xc_\hc(W(C),\,I(C);\,\mc)$). Parity of this cocycle is equal to the parity of $m$.
Finally since the set of all
odd numbers is cofinal in $\mathbb N$ it follows that this
collection determines an element in $H^{m}(\xc_\hc(W(C),\,I(C);\,\mc))$.

In this way we obtain a map
\begin{equation}
\label{eqinvlim}
HC^*_\hc(C,\,\mc)\to H^{*}(\xc_\hc(W(C),\,I(C);\,\mc)).
\end{equation}
In effect one can make one more conjecture
\begin{conj}
Thed map \eqref{eqinvlim} is an isomorphism of the cohomology groups.
\end{conj}
However we shall not need this statement.

\subsection{The cup-products}
Now we are going to define the pairing in Hopf-type cohomology of algebra
$A$ and coalgebra $C$ using the bivariant theories described in this
section. In effect we shall use the maps \eqref{eqbiv1}, \eqref{eqinvlim},
\eqref{eqmapcomp1} (and the corresponding map in cohomology, see Prop. \ref{propmapcomp1})
and the composition product in bivariant
theory. Indeed if $(R,\,I)$ is a quasi-free \hc\/- and C\/- equivariant
extension of $A$ then $HC^p_\hc(A,\,\mc)=
H^{p(\mathrm{mod}\/2)}(\Hom(\xc^p_\hc(R,\,I;\,\mc),\,\Bbbk))=
H^{q}(\xc_\hc(R,\,I;\,\mc),\,\Bbbk)$ (the last group is the group of bivariant cohomology
of tower $\xc_\hc(R,\,I;\,\mc)$
and the constant tower $\Bbbk$). Consider the following composition of maps
\begin{equation}
\label{supermap}
\begin{split}
&HC^p_\hc(C,\,\mc)\otimes HC^q_\hc(A,\,\mc)\stackrel{\eqref{eqinvlim}}{\lar} H^{p}(\xc_\hc(W(C),\,I(C);\,\mc))\otimes HC^q_\hc(A,\,\mc)\\
&\stackrel{\eqref{eqmapcomp1}}{\lar} HC^p(A,\,R,\,I;\hc,\,\mc)\otimes
H^{q}_\hc(A,\,\mc)\stackrel{\eqref{eqbiv1}}{\lar}HC^p(\xc_A,\,\xc_\hc(R,\,I;\,\mc))\otimes
HC_\hc^q(A,\,\mc)\\
&=HC^p(\xc_A,\,\xc_\hc(R,\,I;\,\mc))\otimes
H^{q}(\xc^p_\hc(R,\,I;\,\mc),\,\Bbbk)\to HC^{p+q}(\xc_A,\,\Bbbk)=HC^{p+q}(A).
\end{split}
\end{equation}
The last map in this composition is just the composition product, defined in
\cite{CQ2}. Thus we obtain another cup-product
\begin{equation}
\label{newpair}
HC^p_\hc(C,\,\mc)\otimes
HC^q_\hc(A,\,\mc)\stackrel{\sharp''}{\lar}HC^{p+q}(A).
\end{equation}

It is convenient to write down the cup-product \eqref{newpair} in the terms of
maps $\widetilde\xc^p_\hc\rho_\mc$. We start with the case of even degree elements in $HC^*_\hc(A,\,\mc)$.
For instance, let $[x]\in HC^p(C,\,\mc)$ and $[y]\in HC^{2n}(A,\,\mc)$. Let
$y\in\Hom^0(\xc_\hc^{2n}(RA,\,IA;\,\mc),\,\Bbbk)$ be a cocycle representing
$[y]$ and $x$ be a cocycle in $[x]$. Then as it is mentioned in
dicussion following the proof of the theorem \ref{theoa}, the isomorphism
$\alpha_n$ can be raised to the level of cochains. So we obtain the
following map
\begin{equation}
\label{cupsupeven}
y\sharp''x=y\circ\widetilde\xc^{2n+1}_\hc\rho_\mc(\alpha_n^{-1}(x)):X(B(A))\to
\Bbbk.
\end{equation}
From the properties of all the maps that appear in this formula it
follows that $y\sharp''x$ is a cocycle and that the cohomology class of
this cocycle does not depend on the choices made in its construction (cf. Prop. \ref{propsup}).

In the case of the odd degree elements in $HC^*_\hc(A,\,\mc)$, we can
proceed in the following way. First we observe that the
rows of $T_n\yc$ for the bitower \eqref{diabiv} for $\yc^k_p=\xc_\hc^k(W(C),\,I(C);\,\mc)_p$
can be regarded as the (inverse) towers of the $n$\/-th suspension of the tower, associated
with $\xc_\hc^k(W(C),\,I(C);\mc)$. Here, in analogy with p. 384 of \cite{CQ2} we
put for an inverse tower $\xc=\{\xc_p\}$
\begin{equation}
\label{susp1}
\xc[n]_p^\nu=\xc_{p+n}^{\nu+n},
\end{equation}
with differential $d[n]=(-1)^nd$.
Then the inverse limit of such complexes is equal to the periodic super-complex associated
with $\xc_\hc^k(W(C),\,I(C);\,\mc)[n]$, the $n$\/-th suspension of the
supercomplex of cochain complexes, which is defined by the same formula as
above. Now it is clear that up to a shift of dimension this is the periodic complex of
$\xc_\hc^k(W(C),\,I(C);\mc)$.

On the other hand, since
$\xc_\hc^k(W(C),\,I(C);\,\mc)=X_\hc(W(C);\mc)/F^k_{I(C)}X_\hc(W(C);\,\mc)$, and
the complex $X_\hc(W(C);\,\mc)$ is acyclic, we conclude that instead of the
bitower $\xc_\hc^k(W(C),\,I(C);\,\mc)_p$ one can consider the bitower,
constructed the supercomplexes $F^k_{I(C)}X_\hc(W(C);\,\mc)$. In fact
since the connecting homomorphisms relating the cohomology of
$\xc_\hc^k(W(C),\,I(C);\,\mc)$ and $F^k_{I(C)}X_\hc(W(C);\,\mc)$ commute
with the maps in bitower and establish isomorphisms on the level of
cohomologies, we conclude that the bitower constructed from
$F_{I(C)}X_\hc(W(C);\,\mc)$ has the same cohomology as that of
$\xc_\hc^k(W(C),\,I(C);\,\mc)$.

Thus instead of the elements
$\alpha_n^{-1}(x)\in\xc_\hc^{2n+1}(W(C),\,I(C);\mc)$ one can consider the
corresponding elements $\beta_n^{-1}(x)$ in the $2n+1$\/-st level of
filtration. They fit together to define an element in the limit of bitower
$F_{I(C)}X_\hc(W(C);\,\mc)_p$. On the other hand, this bitower is mapped by
the map $X_\hc(\rho^\sharp,\,\mc)$ constructed in the Proposition
\ref{prneww} to $\Hom(X(B),\,F^{2n}_IX(R,\,\mc))$. Finally we need to
remark that $2n+1$ degree cohomology classes of $A$ correspond to even
cocycles on $F^{2n}_{I(C)}X_\hc(W(C);\,\mc)$. Since all the maps in our
construction of the cup-product $\sharp''$ evidently commute with the
homomorphisms relating the towers $F_{I(C)}X_\hc(W(C);\,\mc)$ and
$\xc_\hc^k(W(C),\,I(C);\,\mc)$ and similar homomorphisms for all the other
towers that appear in the construction of \eqref{supermap} we conclude that
the construction involving $F^*_{I(C)}X_\hc(W(C),\,\mc)$ instead of
$\xc_\hc(W(C),\,I(C);\,\mc)$ yields the same results.

Thus for an odd-degree class $[y]$ in $HC^{2n+1}_\hc(A,\,\mc)$ we can choose
a representing it even cocycle $\tilde y$ in
$F^{2n}_{I(C)}X_\hc(W(C);\,\mc)$; then the class $[x]\sharp''[y]$ is
represented by
\begin{equation}
\label{cupsupodd}
\tilde y\sharp''x= y\circ
X_\hc(\rho^\sharp,\,\mc)(\beta_n^{-1}(x)):X(B(A))\to \Bbbk.
\end{equation}
\begin{prop}
The cup-product $\sharp''$ coincides with $\sharp'$.
\end{prop}
\begin{proof}
From the remark \ref{remmorph} it follows that
the map $\widetilde\xc^p_\hc\rho_\mc$ is an extension of the morphism \eqref{xrho1},
$$
\tilde X^{2n+1}_\hc(\rho^\sharp):\xc_\hc^{2n+1}(W(C),\,I(C);\,\mc)\to\Hom(X(B(A)),\,(\mc\otimes_\hc
R/I^{n+1})_\natural).
$$
Thus from \eqref{cupsupeven} we see that in the case of even degree classes $[y]$,
$[x]\sharp''[y]=[x]\sharp'[y]$.

In the same way the odd-degree case follows from \eqref{cupsupodd} and the simple observation that
the natural projection $X_\hc(R,\,\mc)\to(\mc\otimes_\hc R)_\natural$ sends
$F_I^{2n}X_\hc(R,\,\mc)$ to $I^\hc_n(C,\,\mc)_\natural$ and allows identify
the even cocycles in $F_I^{2n}X_\hc(R,\,\mc)$ with the odd \mc\/-traces
on $A$.
\end{proof}
As a simple consequence of the theory developed in this section one can
derive the following statement
\begin{cor}
There exists a cap-product
\begin{equation}
HC_p(A)\otimes HC^q_\hc(C,\,\mc)\stackrel{\flat}{\lar}HC_{p-q}^\hc(A,\,\mc).
\end{equation}
This product verifies the next equations with respect to $S$\/-operations and
the cup-product $\sharp'$
\begin{align}
\label{eqflat1}
S([z]\flat[x])=S&[z]\flat[x]=[z]\flat S[x],\\
\label{eqflat2}
\langle[z]\flat[x],\,[y]\rangle&=\langle[z],[x]\sharp'[y]\rangle.
\end{align}
Here $\langle,\,\rangle$ is the pairing between the cyclic homology and
cohomology of algebra $A$ (cohomology with coefficients in \mc\ on the left hand side).
\end{cor}
\begin{proof}
The definition of $\flat$ is based on the fact that $HC_*(A)=H_*(\Bbbk,\,\xc_A)$ and
$HC_*^\hc(A,\,\mc)=H_*(\Bbbk,\,\xc_\hc(R,\,I;\,\mc))$ and is similar to that of $\sharp''=\sharp'$,
\eqref{supermap}. In fact one uses the following composition of maps
\begin{equation}
\begin{split}
&HC_p(A)\otimes HC^q_\hc(C,\,\mc)\stackrel{\eqref{eqinvlim}}{\lar}HC_p(A)\otimes H^{q}(\xc_\hc(W(C),\,I(C);\,\mc))\\
&\stackrel{\eqref{eqmapcomp1}}{\lar}HC_p(A)\otimes HC^q(A,\,R,\,I;\hc,\,\mc)
\stackrel{\eqref{eqbiv1}}{\lar}HC_p(A)\otimes HC^q(\xc_A,\,\xc_\hc(R,\,I;\,\mc))\\
&=H_p(k,\,\xc_A)\otimes HC^q(\xc_A,\,\xc_\hc(R,\,I;\,\mc))\to HC^{-p+q}(\xc_A,\,k)=HC_{p-q}(A).
\end{split}
\end{equation}
And the equalities \eqref{eqflat1}, \eqref{eqflat2} follow from the standard
properties of the composition pairing in bivariant cyclic cohomology.
\end{proof}

\section{Other constructions and comparison theorem}

In this section we introduce a new pairing construction in Hopf-cyclic
cohomology with coefficients which generalizes the Quillen's and
Crainic's constructions (see \cite{Quil} and \cite{Crain}) and show
that it coincides with cup-product defined by Khalkhali and
Rangipour in \cite{KhalRan}.



So let $C$ be a coalgebra acting on an algebra $A$, in a way, intertwining the actions of \hc\ on both sides.
Observe that any action of a coalgebra $C$ on an algebra $A$ can be regarded as a map
$$
\rho_C: C\to \Hom(A,A),\quad \rho_C(c)(a) = c(a),\ c\in C, a\in A.
$$
Due to the universal property of Weil algebra $W(C)$, it induces a homomorphism of bigraded bidifferential algebras
$$
\rho_C: W(C) \to \Hom(BA,\Omega),
$$
where coalgebra $BA$ is the bar resolution of $A$ and $\Omega$ is the abbreviation of
$\Omega A$, the
universal differential algebra generated by $A$ (the universal
differential calculus of $A$).
The homomorphism $\rho_C$ intertwines the differential $\delta$ of $BA$ with the
differential $\delta_F$ in $W(C)$ and differential $d$ of $\Omega A$
with $d$ in $W(C)$. If $A$ is an \hc\/-algebra then $\rho_C$ is
an \hc\/-linear map. There's an \hc\/-module structure on $\Hom(BA,\Omega)$,
given by the action of \hc\ on $\Omega A$.

For any SAYD-module \mc\ and an integer $n$ we have a map $\rho_{C,\mc}$
induced by $\rho_C$:
\begin{equation*}
\rho_{C,\mc}: (\mc\otimes_\hc (W(C)^n/\im d))_\natural \to \Hom( BA^\natural,
(\mc\otimes_\hc (\Omega^n/\im d))_\natural).
\end{equation*}
Given a closed graded $(\hc,\mc)$-trace of degree $n$ on $\Omega A$ in the
sense of~\cite{KhalRan}, i.e. a linear functional $\int: (\mc\otimes_\hc (\Omega^n/\im
d))_\natural\to k$, one obtains a chain map of complexes
$$
 \int\circ\rho_{C,\mc}: \left((\mc\otimes_\hc (W(C)^n/\im d))_\natural,\delta_F\right)
\to \left(\Hom(BA^\natural,\Bbbk),\delta\right)
$$
and hence a map in cohomology
\begin{multline}
(\int\circ\rho_{C,\mc})_*: HC^*(C,\mc;n) = H^*((\mc\otimes_\mc (W(C)^n/\im d))_\natural,\delta_F)\\
\to H^*(\Hom(BA^\natural,\Bbbk),\delta) = HC^{*+1}(A).
\end{multline}
Taking into account theorem~\ref{theoa} and~\cite[lemma 3.2]{KhalRan}
that identifies closed $(\hc,\mc)$-traces with Hopf-cyclic cocycles of
$A$, we obtain a pairing:
$$
(\rho_{C,\mc})_*: H^m_H(C,\mc)\otimes H^n_\hc(A,\mc)\to HC^{m+n}(A).
$$

Let us show now that the construction defined above coincides
with cup product of~\cite{KhalRan}.

Denote $W = W(\Bbbk)$ the Weil algebra of the ground field. One always
has a map of bigraded bidifferential algebras
$$ \rho = \rho_k: W \to \Hom(BA,\Omega A).$$
On the other hand, there is a morphism of bigraded bidifferential
algebras
$$
 \sigma : W \to \Hom_\hc(\Omega C, W(C)),
$$
induced by map $\Omega C\stackrel{p}{\rightarrow}\Omega C_0 = C\stackrel{id}{\longrightarrow} C \subset
W(C)$, where $\Omega C$ is universal differential coalgebra of $C$ (see previous sections)
and differentials in $\Hom_\hc(\Omega C, W(C))$ are given by
formulas $df = [d,f]$ and $\delta f = \delta\circ f$.
Finally, thanks to universal property of $\Omega A$ one has a
homomorphism of differential graded algebras
$$
\tau: \Omega A \to \Hom_\hc(\Omega C,\Omega A),
$$
induced by map $\tau:A\to \Hom_\hc(C,A),\  \tau(a)(c) = c(a).$
These homomorphisms give the following diagram (recall that $\Omega$ is an abbreviation for $\Omega A$)
\begin{equation*}
\begin{CD}
W @>\rho>> \Hom(BA,\Omega) @>\tau>> \Hom(BA,\Hom_\hc(\Omega C,\Omega ))
@>>> \Hom_\hc(\Omega C\otimes BA, \Omega)\\
@|        @.                            @.  @|   \\
W @>\sigma>>  \Hom_\hc(\Omega C, W(C)) @>\rho_C>> \Hom_\hc(\Omega C,\Hom(BA,\Omega))
@>>> \Hom_\hc(\Omega C\otimes BA, \Omega).
\end{CD}
\end{equation*}
The diagram is commutative because for the generator $i = i_1$ of
$W$ the upper and the lower row lead to the same map $c\otimes a \mapsto
c(a)$. After factorization we get another diagram of
complexes with differential $\delta$.
$$
\xymatrix{
{ \Hom(BA^\natural,\Omega A_\natural/\im d)}  {\ar[dr]^\tau}\\
& {\hspace{-2cm}\Hom(BA^\natural,\Hom((\mc\otimes_\hc\Omega C)^{\natural, d},(\mc\otimes_\hc\Omega A)_\natural/\im d))}\ar@<-4ex>[d]\\
{W_\natural/\im d}\ar[uu]_\rho \ar[dd]^\sigma & {\hspace{-2cm} \Hom((\mc\otimes_\hc\Omega C)^{\natural, d}\otimes BA^\natural, (\mc\otimes_\hc\Omega A)_\natural/\im d)}\\
& {\hspace{-2cm}\Hom((\mc\otimes_\hc\Omega C)^{\natural, d},\Hom(BA^\natural,(\mc\otimes_\hc\Omega A)_\natural/\im d))}\ar@<4ex>[u]\\
{\Hom((\mc\otimes_\hc\Omega C)^{\natural, d}, (\mc\otimes_\hc W(C))_\natural/\im d).}\ar[ur]_{\rho_{C,\mc}}
}
$$
Here $(\mc\otimes_\hc\Omega C)^{\natural, d}$ is the space of closed
graded $(\hc,\mc)$-cotraces on $\Omega C$ in the sense
of~\cite{KhalRan}, i.e. the set of elements $\xi = \sum_i m_i\otimes \theta_i \in \mc\otimes_\hc \Omega C$
such that $(\id\otimes d)\xi = 0$ and
$$\sum_i (-1)^{|\theta_i^{(1)}||\theta_i^{(2)}|}m_i^{(0)}\otimes \theta_i^{(2)}\otimes m_i^{(-1)}\theta_i^{(1)} =
 \sum_i m_i\otimes \theta_i^{(1)}\otimes \theta_i^{(2)}.$$

Let $\int$ be a closed graded $(\hc,\mc)$-trace of degree $n$  on $\Omega A$
and $\xi$ be a closed graded $(\hc,\mc)$-cotrace of degree $m$ on
$\Omega C$. Evaluating the upper map of the diagram on $\int$ and
$\xi$, we obtain a homomorphism
$$
(\int\cup\xi)\circ\rho: W_\natural/\im d \to \Hom(BA^\natural, \Bbbk),
$$
which is the composition of of $\rho$ and cup product cocycle
$$
\int\cup\xi: \Omega A^{n+m}_\natural/\im d\stackrel{\tau}{\rightarrow}
\Hom((\mc\otimes_\hc\Omega C)_m^{\natural, d},(\mc\otimes_\hc\Omega A)^n_\natural/\im d)
\stackrel{\int\circ \mathbf{ev}_\xi}{\longrightarrow} \Bbbk$$
defined in~\cite{KhalRan}. Take element $cs_{m+n} = iw^{m+n}\in
W$. It is a cocycle in $W_\natural/\im d$ and the map
$$ \rho(cs_{m+n}): BA_{n+m+1}\to \Omega A^{n+m},\quad a_1\otimes
a_2\otimes\dots\otimes a_{m+n+1}\mapsto a_1da_2\dots da_{m+n+1}$$
just identifies closed traces on $\Omega A$ with cyclic cocycles
of $A$. Thus, the cup product of $\int$ and $\xi$ defined in~\cite{KhalRan} can be
expressed as composition cocycle
$$
\int\circ\mathbf{ev}_\xi\circ\tau\circ\rho(cs_{m+n}): BA^\natural_{m+n+1}\to\Bbbk.
$$

On the other hand, the lower row of the diagram shows that the cup
product of $\int$ and $\xi$ coincides with the pairing we defined
of  $\int$ and the element $\sigma(cs_{m+n})(\xi)\in (\mc\otimes_\hc W(C))_\natural/\im
d).$ Remark that due to the obvious inclusion
$$
\Omega_m C = C\otimes (\ker\varepsilon)^m\subset C^{\otimes m+1}\subset
W(C)
$$
$\xi$ can be regarded as element in $(\mc\otimes_\hc W(C))_\natural/\im d$
and in fact it is a cocycle in this complex.
Thus, in order to show that the cup-product of~\cite{KhalRan}
coincides with the  pairing construction of this section it is
sufficient to check the following proposition.

\begin{prop}
$\alpha_{n}[\sigma(cs_{m+n})(\xi)] = \frac {m+1}{m+n+1}[\xi]$ in $HC_\hc^{m}(C,\mc)$
\end{prop}
\begin{proof}
We shall prove this identity by induction on $n$. If $n=0$ then we have
$$
\sigma(cs_m)(c_0\otimes c_1\otimes\dots\otimes c_m) = c_0c_1\dots c_m
$$
for every $c_0\otimes c_1\otimes\dots\otimes c_m\in\Omega_m C$. So
the identity is true. In order to establish the step of induction
consider the following diagram
$$
\xymatrix{
{0} \ar[r]& {\qquad W^{n+1}/b,d \quad} \ar[d]^{\mathbf{ev}_\xi} \ar[r]^b &{\qquad W^n/d,(1-\kappa)\quad} \ar[d]^{\mathbf{ev}_\xi} \ar[r]&{\qquad\qquad\quad } \\
{0} \ar[r] & {(\mc\otimes_\hc W(C)^{n+1})/b,d}\ar[r]^b& {(\mc\otimes_\hc W(C)^n)/d,(1-\kappa)}\ar[r] &{\qquad\qquad\quad}\\
&         &                                    {\quad\qquad\qquad} \ar[r]&   {\qquad W^n/b,d \quad} \ar[r] \ar[d]^{\mathbf{ev}_\xi}& 0\\
&           &                                  {\quad\qquad\qquad}\ar[r]&      {(\mc\otimes_\hc W(C)^n)/b,d} \ar[r] &0.
}
$$
It induces a commutative diagram for boundary maps:
\begin{equation*}
\begin{CD}
H^*(W^n_\natural/\im d, \delta) @>\varphi>> H^*(W^{n+1}_\natural/\im d,
\delta)\\
@VV\mathbf{ev}_\xi V @VV\mathbf{ev}_\xi V\\
H^*_\hc(C,\mc; n)  @>\varphi>> H^*_\hc(C,\mc; n+1)
\end{CD}
\end{equation*}
Show now that $\varphi[cs_n] = \frac{n+1}{n+2}[cs_{n+1}]\in H^*(W^{n+1}_\natural/\im d,
\delta)$.
Indeed, we have the equality in $W/\im(1-\kappa)$
\begin{multline*}
d(\sum_{k=0}^{n-1}i^2w^kiw^{n-k-1}) = \sum_{k=0}^{n-1}
(wiw^kiw^{n-k-1} -iw^{k+1}iw^{n-k-1}+i^2w^n) =
(n+1)i^2w^n-iw^ni,
\end{multline*}
that implies $iw^ni = (n+1)i^2w^n$ in $W/\im d+\im(1-\kappa)$ and
\begin{eqnarray*}
b(iw^{n+1}) = -iw^ni-i^2w^n = -(n+2)i^2w^n,\\
\delta(iw^n) = -iw^ni = -(n+1)i^2w^n.
\end{eqnarray*}
Hence, the equality holds $\varphi[\sigma(cs_{n+1})(\xi)] =
\frac{n+1}{n+2}[\sigma(cs_n)(\xi)]$ and
\begin{multline*}
\alpha_{n+1}[\sigma(cs_{m+n+1})(\xi)] = \alpha_n\circ\varphi[\sigma(cs_{m+n+1})(\xi)] =
\frac{m+n+1}{m+n+2}\alpha_n[\sigma(cs_{n+m})(\xi)] \\
=\frac{m+n+1}{m+n+2}\cdot\frac {m+1}{m+n+1}[\xi] = \frac {m+1}{m+n+2}[\xi].
\end{multline*}

\end{proof}


\begin{thebibliography}{99}
\bibitem{HKRS}{P.M.Hajac, M.Khalkali, B.Rangipour, M.Sommerh\"auser: {\sc Hopf-cyclic homology and cohomology with coefficients}; C. R. Math. Acad. Sci. Paris {\bf 338} (2004), no. 9, 667-672 (also available as preprint
arXiv:math.KT/0306288 v.2)}
\bibitem{CQ2}{J.Cuntz, D.Quillen: {\sc Cyclic homology and
nonsingularity}; J. Amer. Math. Soc. {\bf 8} n.2 (1995) 373-442}
\bibitem{CQ1}{J.Cuntz, D.Quillen: {\sc Algebra extensions and
nonsingularity}; J. Amer. Math. Soc. {\bf 8} n.2 (1995) 251-289}
\bibitem{CQ3}{J.Cuntz, D.Quillen: {\sc Excision in bivariant  periodic cyclic
cohomology}; Invent. Math. {\bf 127} (1987) 67-98}
\bibitem{Good}{T.G.Goodwillie:{\sc Cyclic homology, derivations,
and the free loopspace}; Topology {\bf 24} n. 2 (1985) 187--215}
\bibitem{JK}{J.D.S.Jones, C.Kassel {\sc Bivariant cyclic theory}; $K$\/-theory,
{\bf 3} (1989), 339-365}
\bibitem{Crain}{M.Crainic: {\sc Cyclic cohomology of Hopf algebras}; J. Pure Appl. Algebra {\bf 166} (2002) 29-66}
\bibitem{CM1}{A.Connes, H.Moscovici: {\sc Hopf algebras, cyclic cohomology and the transverse index theorem}; Commun. Math. Phys. {\bf 198} (1998), 199-246}
\bibitem{CM2}{A.Connes, H.Moscovici: {\sc Cyclic cohomology and Hopf algebras}; Lett. Math. Phys. {\bf 52} (2000) 97-108}
\bibitem{Teill}{R.Taileffer: {\sc Cyclic homology of Hopf algebras}; K-Theory {\bf 24} (2001) 69-85}
\bibitem{Iran1}{M.Akbarpour, M.Khalkhali: {\sc Equivariant cyclic cohomology of Hopf module algebras}; J. reine angew. Math. {\bf 559} (2003) 137-152}
\bibitem{Iran2}{M.Khalkhali, B.Rangipour: {\sc Invariant cyclic homology}; preprint 
arXiv:math.KT/0207118}
\bibitem{Iran3}{M.Khalkhali, B.Rangipour: {\sc A new cyclic module for Hopf algebras}; K-theory {\bf 27} (2002) 111-131}
\bibitem{Loday}{J.L.Loday: {\sc Cyclic homology}; Springer-Verlag, 1998}
\bibitem{KhalRan}{M.Khalkhali, B.Rangipour: {\sc Cup Products in Hopf-Cyclic Cohomology}; available as preprint at arXiv:math.QA/0411003 v1}
\bibitem{KayKhal}{A.Kaygun, M.Khalkhali: {\sc Excision in Hopf cyclic homology}; available as preprint at arXiv:math.KT/0511026 v1}
\bibitem{Karoubi}{M. Karoubi: {\sc Cohomologie cyclique et K-th\'eorie}; Ast\'erisque, {\bf 149}, 1987}
\bibitem{Quil}{D.Quillen: {\sc Chern-Simons forms and cyclic cohomology}; The interface of Mathematics and particle Physics, (Oxford, 1988) 117-134}
\bibitem{Quil2}{D.Quillen: {\sc Algebra cochains and cyclic cohomology}; Publ. Math. IHES, {\bf 68}, (1989) 139-174}
\bibitem{Shar}{G.Sharygin: {\sc A new construction of characteristic classes for noncommutative algebraic principal bundles}; Banach Center Publ. {\bf 61} (2003) 219-230}
\bibitem{nashe}{I.Nikonov, G.Sharygin: {\sc On the Hopf-type cyclic cohomology with coefficients}; submitted for publication in the proceedings of
the conference ``$C^*$\/-algebras and Elliptic Theory, I", B\c{e}dlewo, 2004}
\end{thebibliography}
\end{document}